\documentclass[10pt]{article}
\usepackage[latin9]{inputenc}
\usepackage[a4paper]{geometry}
\usepackage{color}
\usepackage{float}
\usepackage{mathrsfs}
\usepackage{amsmath}
\usepackage{microtype}

\makeatletter

\providecommand{\tabularnewline}{\\}
\floatstyle{ruled}
\newfloat{algorithm}{tbp}{loa}
\providecommand{\algorithmname}{Algorithm}
\floatname{algorithm}{\protect\algorithmname}

\usepackage{amssymb}
\usepackage{algorithm,algpseudocode}
\usepackage{pgf}
\graphicspath{{plots/}}

\makeatother

\usepackage{listings}
\begin{document}
\title{A low-rank complexity reduction algorithm for the high-dimensional
kinetic chemical master equation}
\author{Lukas Einkemmer\thanks{Department of Mathematics, Universität Innsbruck, Innsbruck, Tyrol, Austria}\;\,\footnote{lukas.einkemmer@uibk.ac.at} \and Julian Mangott\footnotemark[1] \and Martina Prugger\thanks{Department of Biochemistry, Universität Innsbruck, Innsbruck, Tyrol, Austria}
}
\maketitle
\begin{abstract}
It is increasingly realized that taking stochastic effects into account
is important in order to study biological cells. However, the corresponding
mathematical formulation, the chemical master equation (CME), suffers
from the curse of dimensionality and thus solving it directly is not
feasible for most realistic problems. In this paper we propose a dynamical
low-rank algorithm for the CME that reduces the dimensionality of
the problem by dividing the reaction network into partitions. Only
reactions that cross partitions are subject to an approximation error
(everything else is computed exactly). This approach, compared to
the commonly used stochastic simulation algorithm (SSA, a Monte Carlo
method), has the advantage that it is completely noise-free. This
is particularly important if one is interested in resolving the tails
of the probability distribution. We show that in some cases (e.g.~for
the lambda phage) the proposed method can drastically reduce memory
consumption and run time and provide better accuracy than SSA.
\end{abstract}

\section{Introduction}

Chemical kinetics is an indispensable tool in order to understand
reaction networks that govern, for example, the chemical processes
inside a biological cell. The fundamental mathematical description
of such systems is the chemical master equation (CME). However, since
each chemical species adds a dimension to the CME, solving it numerically
is extremely expensive. More precisely, the memory required and the
computational cost scales exponentially in the number of species.
This is often referred to as the \emph{curse of dimensionality}. As
a consequence, reduced models that only take averaged population numbers
into account are most commonly used \cite{Chen_2010}. This assumption
results in a set of ordinary differential equations (ODE) that can
then be solved at low computational cost. ODE models are also called
\emph{deterministic}, owing to the fact that they only give averaged
values and thus neglect both the inherent stochasticity of the system
as well as the discrete nature of population numbers. It is increasingly
realized, however, that both are required in order to describe many
important features in biological systems \cite{Tonn_2019,Grima_2008,Niepel_2009,Paszek_2010}.
Thus computing a solution of the full chemical master equation is
required in order to understand such systems.

Directly solving the chemical master equation for realistic system
sizes is either very costly or prohibitive in terms of both memory
and computational cost (primarily due to the curse of dimensionality).
The most commonly used approach currently is the stochastic simulation
algorithm (SSA; see, e.g., \cite{Gillespie_1976,Harris_2006}). The
SSA is a Monte Carlo approach that simulates individual trajectories
of the system. While one such sample, owing to the inherent randomness,
does not tell us much useful information, repeating it many times
allows us to collect a statistic of the most likely outcomes of the
system. As a Monte Carlo method SSA does not suffer from the curse
of dimensionality. However, it only converges slowly (as $1/\sqrt{N}$,
where $N$ is the number of samples) and is very noisy if not enough
samples are used. The latter is a phenomenon where even if the probability
density function is perfectly smooth, the algorithm approximates it
by a jagged line. This, in particular, is an issue for the tail of
the distribution, where the noise can completely bury the physical
behavior of the system.

In this paper we propose a method that directly reduces the dimensionality
of the problem by using a low-rank approximation. In this approach,
lower-dimensional basis functions (which require far less memory to
store) are combined in order to obtain an approximation to the high-dimensional
problem. For the degrees of freedom in the low-rank approximation
(also called the \emph{low-rank factors}) we derive evolution equations
that are then used to advance the approximation forward in time. This
dynamical low-rank approach dates back to early work in quantum mechanics
(see, e.g., \cite{Meyer_1990,Meyer_2009,Lubich_2008}) and a number
of important mathematical advances in constructing and analyzing such
methods have been made more recently \cite{Lubich_2014,Kusch_2023,Ceruti_2022b,Ceruti_2022a,Einkemmer_2022,Einkemmer_2022a,Ceruti_2022,Ding_2021,Einkemmer_2021b}.
In the quantum mechanics context usually single-orbital basis functions,
which only depend on the coordinates of a single electron, are combined
to obtain an approximation to the high-dimensional wave function.
In \cite{Jahnke_2008} this idea has been directly applied to the
chemical master equation. The problem with that approach, however,
is that each of the low-rank factors are only allowed to depend on
a single species. It is doubtful that in biological applications,
given the intricate structures of complex biological networks \cite{Barabasi_2009},
we can consider each species independently, while still obtaining
an accurate approximation with a small rank.

What we propose in this paper is to divide a reaction network into
two partitions. The low-rank factors are then allowed to depend on
all species in their respective partition. Thus, all reactions inside
of a partition are treated exactly. An approximation is only performed
if a reaction crosses the partition boundary. This allows us to keep
species that tightly couple to each other together without introducing
any error, while still taking advantage of the computational and memory
savings of the dynamical low-rank approach. We emphasize that computational
savings are not only due to lower-dimensional low-rank factors, but
also depend crucially on how small the rank (i.e.~the number of such
low-rank factors used) can be chosen while still maintaining accurate
results. A similar approach has been used in \cite{Prugger_2023}
for Boolean models in biology. In the present work we extend this
to the full kinetic chemical master equation. Let us also note that
similar approaches have been used for problems in plasma physics (see,
e.g., \cite{Einkemmer_2018,Cassini_2021,Coughlin_2022,Einkemmer_2020})
and radiatiation transport (see, e.g., \cite{Peng_2020,Peng_2021,Kusch_2021,Einkemmer_2021,Einkemmer_2021a,Kusch_2022}).
In this case the partitioning is also based on the underlying physical
problem (either a decomposition into spatial and velocity scales,
as in \cite{Peng_2020,Kusch_2021,Einkemmer_2020,Einkemmer_2018},
or in coordinates parallel and perpendicular to the magnetic field,
as in \cite{Einkemmer_2023}). Our view is that in biological applications
there are a multitude of different reaction networks and, in general,
for each a different partitioning will give optimal results.

The remainder of the paper is structured as follows. In section \ref{sec:cme}
we introduce the chemical master equation and set our notation. The
dynamical low-rank approximation is then described in detail in section
\ref{sec:DLR-approximation}. In section \ref{subsec:implementation-general-remarks}
we discuss the steps that are necessary in order to obtain an efficient
implementation. In section \ref{sec:numerical-experiments} we investigate
the accuracy and efficiency of the method for a number of examples.
In particular, we show that for a lambda phage model the proposed
algorithm is more accurate than SSA and drastically reduces the required
run time. Finally, we conclude in section \ref{sec:Outlook}.

\section{Chemical master equation\label{sec:cme}}

A well-stirred chemical reaction system of $N$ species $S_{1},\ldots,S_{N}$
is interacting through $M$ reaction channels $R_{1},\ldots,R_{M}$.
In the stochastic description the system is represented by a random
variable $\mathcal{X}(t)=\left(\mathcal{X}_{1}(t),\ldots,\mathcal{X}_{N}(t)\right)$
on the discrete state space $\mathbb{N}_{0}^{N}$, where the entries
$\mathcal{X}_{i}(t)$ denote the population number (i.e. number of
molecules) of the $i$-th species at time $t$. The probability density

\[
P(t,x)=\mathbb{P}(\mathcal{X}_{1}(t)=x_{1},\ldots,\mathcal{X}_{N}(t)=x_{N}),\quad x=(x_{1},\ldots,x_{N})\in\mathbb{N}_{0}^{N},
\]
where $x$ are the population numbers, is the solution of the kinetic
chemical master equation (CME)

\begin{equation}
\partial_{t}P(t,x)=\sum_{\mu=1}^{M}\left(a_{\mu}(x-\nu_{\mu})P(t,x-\nu_{\mu})-a_{\mu}(x)P(t,x)\right).\label{eq:CME}
\end{equation}

By defining the linear operator
\begin{equation}
\left(\mathcal{\mathscr{A}}P(t,\cdot)\right)(x)=\sum_{\mu=1}^{M}\left(a_{\mu}(x-\nu_{\mu})P(t,x-\nu_{\mu})-a_{\mu}(x)P(t,x)\right),\label{eq:linear_operator}
\end{equation}
the CME can be concisely written as 
\[
\partial_{t}P(t,\cdot)=\mathcal{\mathscr{A}}P(t,\cdot).
\]

The stoichiometric vector $\nu_{\mu}=(\nu_{\mu,1},\ldots,\nu_{\mu,N})\in\mathbb{Z}^{N}$
describes the population change caused by reaction $\mu$. The propensity
function $a_{\mu}(x):\mathbb{N}_{0}^{N}\to[0,\,\infty)$ for reaction
channel $R_{\mu}$ can be interpreted as a transition probability
$T_{\mu}(x+\nu_{\mu}\mid x).$ Note that the arguments $x-\nu_{\mu}$
in the first term of the right-hand side of equation~(\ref{eq:CME})
have to be omitted when they become negative (there can be no physical
reaction that reduces the population number to negative values). The
term ``kinetic'' indicates that the population number can be of
any natural number (including 0), in contrast to models that treat
only boolean states (where a species can be either ``activated''
or ``not activated''; see, e.g., \cite{Clarke_2020,Zanudo_2018,Yachie-Kinoshita_2018,Prugger_2023}).
For more details on the CME in general, we refer the reader to, e.g.,
\cite{Gillespie_1976,Gillespie_1992,Gardiner_2004}.

We will illustrate these concepts with a simple example. The bimolecular
reaction $A+B\rightleftarrows C$ has propensity functions $a_{f}(x)$
for the forward and $a_{b}(x)$ for the backward reaction, with $x=(x_{A},x_{B},x_{C})$.
The propensity functions describe how likely the associated reaction
occurs for the given population number. In the forward reaction a
particle of species $A$ reacts with a particle of species $B$ and
yields one of species $C$, so $\nu_{f}=(-1,-1,1)$. Similarly, the
stoichiometric vector for the backward reaction is $\nu_{b}=(1,1,-1)$.
The entire CME therefore reads as 
\[
\partial_{t}P(t,x)=a_{f}(x-\nu_{f})P(t,x-\nu_{f})+a_{b}(x-\nu_{b})P(t,x-\nu_{b})-\left(a_{f}(x)+a_{b}(x)\right)P(t,x).
\]

An important physical relation is the conservation of probability
(also called \emph{conservation of mass}), 
\begin{equation}
\sum_{x\in\mathbb{N}_{0}^{N}}P(t,x)=\sum_{x\in\mathbb{N}_{0}^{N}}P(0,x)=1,\label{eq:mass}
\end{equation}
 which can be directly derived from the CME. To see this, we integrate
equation~(\ref{eq:CME}) over time and perform a summation over $x$,
yielding 
\begin{equation}
\sum_{x\in\mathbb{N}_{0}^{N}}\left(P(t,x)-P(0,x)\right)=\text{\ensuremath{\int_{0}^{t}}}\sum_{\mu=1}^{M}\sum_{x\in\mathbb{N}_{0}^{N}}\left(a_{\mu}(x-\nu_{\mu})P(\tilde{t},x-\nu_{\mu})-a_{\mu}(x)P(\tilde{t},x)\right)\,\mathrm{d}\tilde{t}.\label{eq:mass_int}
\end{equation}
The right-hand-side of this equation vanishes since
\begin{equation}
\sum_{\mu=1}^{M}\sum_{x\in\mathbb{N}_{0}^{N}}a_{\mu}(x-\nu_{\mu})P(t,x-\nu_{\mu})=\sum_{\mu=1}^{M}\sum_{x\in\mathbb{N}_{0}^{N}}a_{\mu}(x)P(t,x),\label{eq:substitution}
\end{equation}
which implies the desired result.

\section{Dynamical low-rank approximation\label{sec:DLR-approximation}}

Solving the full CME is not possible in most cases due to the curse
of dimensionality. Even the memory requirement for storing the full
probability density function with a finite number $\tilde{n}$ of
possible population numbers for the $N$ species scales with $\mathcal{O}(\tilde{n}^{N})$.
Therefore, we have to reduce the system size in order to solve the
CME using currently available hardware. We will do this using a dynamical
low-rank approximation. The main idea is to split the species (and
thus the reaction network) into two partitions. Reaction pathways
lying within a partition are treated exactly, while reaction pathways
that cross the two partitions are taken into account in an approximate
way.

More specifically, we separate the reaction network into two partitions,
such that there are $m_{1}<N$ species lying in partition 1. We write
the population numbers as $x=(x_{(1)},x_{(2)}),$ such that $x_{(1)}=(x_{1},\ldots,x_{m_{1}})$
and $x_{(2)}=(x_{m_{1}+1},\ldots,x_{N})$ are the population numbers
in partition 1 and 2, respectively. In the following, we will denote
all tuples belonging to the first or second partition by parenthesized
indices, i.e. by $(1)$ or $(2)$. Then the CME reads as
\begin{align}
\partial_{t}P(t,x_{(1)},x_{(2)}) & =\sum_{\mu=1}^{M}a_{\mu}(x_{(1)}-\nu_{\mu,(1)},x_{(2)}-\nu_{\mu,(2)})P(t,x_{(1)}-\nu_{\mu,(1)},x_{(2)}-\nu_{\mu,(2)})\label{eq:CME_2p}\\
 & \qquad-\sum_{\mu=1}^{M}a_{\mu}(x_{(1)},x_{(2)})P(t,x_{(1)},x_{(2)}).\nonumber 
\end{align}
The dynamical low-rank (DLR) approximation of $P$ is given by
\begin{equation}
P(t,x_{(1)},x_{(2)})\approx\sum_{i,j=1}^{r}X_{i}^{1}(t,x_{(1)})S_{ij}(t)X_{j}^{2}(t,x_{(2)}),\label{eq:dlra}
\end{equation}
where $S_{ij}\in\mathbb{R}$ is the coefficient matrix and $r$ is
called the \emph{rank of the representation}. The dependency of $P$
on $x_{(1)}$ and $x_{(2)}$ is approximated by the basis functions
$\{X_{i}^{1}:\,i=1,\ldots,r\}$ and $\{X_{j}^{2}:\,j=1,\ldots,r\}$.
These functions depend on time $t$ but only on the population numbers
in partition 1, $x_{(1)}\in\mathbb{N}_{0}^{m_{1}}$ or on the population
numbers in partition 2, $x_{(2)}\in\mathbb{N}_{0}^{m_{2}}$ (with
$m_{2}=N-m_{1}$). The crucial benefit of this approach is that the
memory requirements for storing the low-rank factors $X_{i}^{1}(t,x_{(1)})$,
$S_{ij}(t)$ and $X_{j}^{2}(t,x_{(2)})$ scales with $\mathcal{O}((\tilde{n}^{m_{1}}+\tilde{n}^{m_{2}})\cdot r+r^{2})$.
As $r$ is usually small, the memory requirements are reduced drastically
compared to the full probability density, which would require $\mathcal{O}(\tilde{n}^{m_{1}+m_{2}})$.
In the following we will add additional constraints so as to obtain
uniqueness of the representation given by equation (\ref{eq:dlra})
and derive an algorithm for computing $X_{i}^{1}(t,x_{(1)})$, $S_{ij}(t)$
and $X_{j}^{2}(t,x_{(2)})$.

Let us assume that $S$ is invertible and $X_{i}^{1}$ and $X_{i}^{2}$
obey the orthogonality and gauge conditions

\begin{align}
\langle X_{i}^{1},X_{j}^{1}\rangle_{1}=\delta_{ij} & \quad\textrm{and}\quad\langle X_{i}^{2},X_{j}^{2}\rangle_{2}=\delta_{ij},\qquad\textrm{(orthogonality)},\label{eq:ortho}\\
\langle X_{i}^{1},\partial_{t}X_{j}^{1}\rangle_{2}=0 & \quad\textrm{and}\quad\langle X_{i}^{1},\partial_{t}X_{j}^{1}\rangle_{2}=0,\qquad\textrm{(gauge condition)},\label{eq:gauge}
\end{align}
where $\delta_{ij}$ denotes the Kronecker delta and $\langle\cdot,\cdot\rangle_{k}$
the inner product on $\ell^{2}(\mathbb{N}_{0}^{m_{k}})$ $(k=1,2)$.
Then the approximation $P\in\ell^{2}(\mathbb{N}_{0}^{N})$ is unique
(see, e.g., \cite{Koch_2007,Einkemmer_2018}) and lies for any time
$t$ in the low-rank manifold 
\begin{align*}
\mathcal{M}= & \bigg\{ P\in\ell^{2}(\mathbb{N}_{0}^{N}):\,P(x_{(1)},x_{(2)})=\sum_{i,j=1}^{r}X_{i}^{1}(x_{(1)})S_{ij}X_{j}^{2}(x_{(2)}),\\
 & \quad\textrm{with invertible}\:S=S_{ij}\in\mathbb{R}^{r\times r},X_{i}^{k}\in\ell^{2}(\mathbb{N}_{0}^{m_{k}})\:\textrm{and}\:\langle X_{i}^{k},X_{j}^{k}\rangle_{k}=\delta_{ij}\;(k=1,2)\bigg\}
\end{align*}
with tangent space
\begin{align*}
\mathcal{T}_{P}\mathcal{M}= & \bigg\{\dot{P}\in\ell^{2}(\mathbb{N}_{0}^{N}):\,\dot{P}(x_{(1)},x_{(2)})=\sum_{i,j=1}^{r}\left(\dot{X}_{i}^{1}(x_{(1)})S_{ij}X_{j}^{2}(x_{(2)})+X_{i}^{1}(x_{(1)})\dot{S}_{ij}X_{j}^{2}(x_{(2)})+X_{i}^{1}(x_{(1)})S_{ij}\dot{X}_{j}^{2}(x_{(2)})\right),\\
 & \quad\textrm{with}\:\dot{S}\in\mathbb{R}^{r\times r},\dot{X}_{i}^{k}\in\ell^{2}(\mathbb{N}_{0}^{m_{k}})\:\textrm{and}\:\langle X_{i}^{k},\dot{X}_{j}^{k}\rangle_{k}=0\;(k=1,2)\bigg\},
\end{align*}
\textcolor{black}{where dotted quantities denote the formal derivative
with respect to time.} Using the orthogonality (\ref{eq:ortho}),
the gauge conditions (\ref{eq:gauge}) and the linear operator defined
in equation (\ref{eq:linear_operator}), we then obtain the following
relations
\begin{equation}
\begin{aligned}\partial_{t}S_{ij} & =\left\langle X_{i}^{1}X_{j}^{2},\mathcal{\mathscr{A}}\sum_{i,j=1}^{r}X_{i}^{1}S_{ij}X_{j}^{2}\right\rangle _{1,2},\\
\sum_{j=1}^{r}S_{ij}\partial_{t}X_{j}^{2} & =\left\langle X_{i}^{1},\mathcal{\mathscr{A}}\sum_{i,j=1}^{r}X_{i}^{1}S_{ij}X_{j}^{2}\right\rangle _{1}-\sum_{j=1}^{r}\partial_{t}S_{ij}X_{j}^{2},\\
\sum_{i=1}^{r}S_{ij}\partial_{t}X_{i}^{1} & =\left\langle X_{j}^{2},\mathcal{\mathscr{A}}\sum_{i,j=1}^{r}X_{i}^{1}S_{ij}X_{j}^{2}\right\rangle _{2}-\sum_{i=1}^{r}X_{i}^{1}\partial_{t}S_{ij}.
\end{aligned}
\label{eq:projections}
\end{equation}

In principle we can solve this set of equations and thus we can determine
the time evolution of the low-rank factors. However, if, e.g.~a classic
Runge--Kutta method is applied to equation (\ref{eq:projections}),
we need to invert $S$. If $S$ has small singular values, inverting
it is numerically very ill-conditioned. If, on the other hand, $S$
has only large singular values, then the approximation is very inaccurate
(this corresponds to the case where the rank $r$ has been chosen
too small to obtain an accurate approximation). This has been realized
early in the development of such schemes, with regularization being
a somewhat unsatisfactory remedy (see, e.g., \cite{Lubich_2008,Meyer_2009}).
In the seminal paper \cite{Lubich_2014} a projector splitting scheme
was introduced that avoids the inversion of $S$ and thus results
in a method that is robust with respect to the presence of small singular
values. Later the basis updating Galerkin (BUG, also called the unconventional
integrator) approach was introduced in \cite{Ceruti_2022} and improved
in \cite{Ceruti_2022a}. Any of these robust integrators would be
suitable for the task at hand. However, since \cite{Prugger_2023}
(for the Boolean case) and \cite{Einkemmer_2023} (for a kinetic problem
from plasma physics) seems to indicate that for reversible problems
the projector splitting integrator seems to be more accurate, consumes
less memory and incurs less computational cost, we will mostly focus
on this approach here. Let us, however, duly note that the integrators
are very similar in the sense that the building blocks we derive below
can also be used easily to implement any of the variants of the BUG
integrator.

We can write equation (\ref{eq:projections}) as
\begin{equation}
\partial_{t}P=\mathscr{P}(P)\mathcal{\mathscr{A}}\sum_{i,j=1}^{r}X_{i}^{1}S_{ij}X_{j}^{2},\label{eq:projector-equation}
\end{equation}
where $P$ is given by the low-rank approximation in equation (\ref{eq:dlra})
and $\mathscr{P}(P)$ is the projector onto the tangent space $\mathcal{T}_{P}\mathcal{M}$
\[
\mathscr{P}(P)g=\sum_{j=1}^{r}\langle X_{j}^{2},g\rangle_{2}X_{j}^{2}-\sum_{i,j=1}^{r}X_{i}^{1}\langle X_{i}^{1}X_{j}^{2},g\rangle_{1,2}X_{j}^{2}+\sum_{i=1}^{r}X_{i}^{1}\langle X_{i}^{1},g\rangle_{1},
\]

for more details see, e.g., \cite{Lubich_2014,Einkemmer_2018}. The
idea of the projector splitting integrator is to treat each term in
the projector separately. That is, we split equation (\ref{eq:projector-equation})
into the following three parts
\begin{align}
\partial_{t}P & =\sum_{j=1}^{r}\left\langle X_{j}^{2},\mathcal{\mathscr{A}}\sum_{i,j=1}^{r}X_{i}^{1}S_{ij}X_{j}^{2}\right\rangle _{2}X_{j}^{2},\label{eq:LT_1}\\
\partial_{t}P & =-\sum_{i,j=1}^{r}X_{i}^{1}\left\langle X_{i}^{1}X_{j}^{2},\mathcal{\mathscr{A}}\sum_{i,j=1}^{r}X_{i}^{1}S_{ij}X_{j}^{2}\right\rangle _{1,2}X_{j}^{2},\label{eq:LT_2}\\
\partial_{t}P & =\sum_{i=1}^{r}X_{i}^{1}\left\langle X_{i}^{1},\mathcal{\mathscr{A}}\sum_{i,j=1}^{r}X_{i}^{1}S_{ij}X_{j}^{2}\right\rangle _{1}.\label{eq:LT_3}
\end{align}

We will now explain the algorithm for computing $X_{i}^{1}(t,x_{(1)})$,
$S_{ij}(t)$ and $X_{j}^{2}(t,x_{(2)})$ by using the first-order
Lie-Trotter splitting. The initial value for the algorithm is given
by 
\[
P(0,x_{(1)},x_{(2)})=\sum_{i,j=1}^{r}X_{0,i}^{1}(x_{(1)})S_{0,ij}X_{0,j}^{2}(x_{(2)}).
\]

In the first step of the algorithm, we solve equation~(\ref{eq:LT_1}).
We write 
\[
P(t,x_{(1)},x_{(2)})=\sum_{j=1}^{r}K_{j}(t,x_{(1)})X_{j}^{2}(t,x_{(2)}),\quad\textrm{with}\quad K_{j}(t,x_{(1)})=\sum_{i=1}^{r}X_{i}^{1}(t,x_{(1)})S_{ij}(t),
\]
this step is therefore commonly called the $K$ \emph{step}. Inserting
this expression into equation~(\ref{eq:LT_1}) yields an equation
whose solution is given by the \emph{time-independent} functions $X_{j}^{2}(t,x_{(2)})=X_{j}^{2}(0,x_{(2)})=X_{0,j}^{2}(x_{(2)})$
(see, e.g., \cite{Einkemmer_2018,Lubich_2014}). After applying an
inner product $\langle X_{i}^{2}(x_{(2)}),\cdot\rangle_{2}$ and using
the orthogonality condition (\ref{eq:ortho}), we further obtain
\begin{equation}
\partial_{t}K_{i}(t,x_{(1)})=\sum_{\mu=1}^{M}\sum_{j=1}^{r}\left(c_{ij}^{1,\mu}(x_{(1)})K_{j}(t,x_{(1)}-\nu_{\mu,(1)})-d_{ij}^{1,\mu}(x_{(1)})K_{j}(t,x_{(1)})\right),\label{eq:K}
\end{equation}
with the \emph{time-independent} coefficients 
\begin{equation}
\begin{aligned}c_{ij}^{1,\mu}(x_{(1)}) & =\langle X_{0,i}^{2}(x_{(2)}),a_{\mu}(x_{(1)}-\nu_{\mu,(1)},x_{(2)}-\nu_{\mu,(2)})X_{0,j}^{2}(x_{(2)}-\nu_{\mu,(2)})\rangle_{2},\\
d_{ij}^{1,\mu}(x_{(1)}) & =\langle X_{0,i}^{2}(x_{(2)}),a_{\mu}(x_{(1)},x_{(2)})X_{0,j}^{2}(x_{(2)})\rangle_{2}.
\end{aligned}
\label{eq:K_coeff}
\end{equation}
The coefficients can be simplified when the propensity function $a_{\mu}(x_{(1)},x_{(2)})$
factorizes in its arguments or depends only on a subset of the population
numbers, as we will investigate later. We now integrate equation~(\ref{eq:K})
with the initial value 
\[
K_{j}(0,x_{(1)})=\sum_{i=1}^{r}X_{0,i}^{1}(x_{(1)})S_{0,ij}
\]
until time $\tau$ to obtain $K_{1,j}(x_{(1)})=K_{j}(\tau,x_{(1)})$.
Then, we perform a QR decomposition 
\[
K_{1,j}(x_{(1)})=\sum_{i=1}^{r}X_{1,i}^{1}(x_{(1)})\hat{S}_{ij},
\]
which gives orthonormal functions $X_{1,i}^{1}$ (remember that on
the low-rank manifold the orthogonality condition (\ref{eq:ortho})
has to be fulfilled) and the matrix $\hat{S}_{ij}$.

In the second step of the algorithm we proceed in a similar way for
equation~(\ref{eq:LT_2}) and notice that the solution is given by
time-independent functions $X_{i}^{1}(t,x_{(1)})=X_{i}^{1}(\tau,x_{(1)})=X_{1,i}^{1}(x_{(1)})$
and $X_{j}^{2}(t,x_{(2)})=X_{j}^{2}(0,x_{(2)})=X_{0,j}^{2}(x_{(2)})$.
After a similar calculation as for the $K$ step, we obtain the central
equation for the $S$ \emph{step}
\begin{equation}
\partial_{t}S_{ij}(t)=-\sum_{k,l=1}^{r}S_{kl}(t)\left(e_{ijkl}-f_{ijkl}\right),\label{eq:S}
\end{equation}
with the time-independent coefficients 
\begin{equation}
\begin{aligned}e_{ijkl} & =\sum_{\mu=1}^{M}\langle X_{1,i}^{1}(x_{(1)})X_{0,j}^{2}(x_{(2)}),a_{\mu}(x_{(1)}-\nu_{\mu,(1)},x_{(2)}-\nu_{\mu,(2)})X_{1,k}^{1}(x_{(1)}-\nu_{\mu,(1)})X_{0,l}^{2}(x_{(2)}-\nu_{\mu,(2)})\rangle_{1,2},\\
f_{ijkl} & =\sum_{\mu=1}^{M}\langle X_{1,i}^{1}(x_{(1)})X_{0,j}^{2}(x_{(2)}),a_{\mu}(x_{(1)},x_{(2)})X_{1,k}^{1}(x_{(1)})X_{0,l}^{2}(x_{(2)})\rangle_{1,2}.
\end{aligned}
\label{eq:S_coeff}
\end{equation}
Note the minus sign in front of the right-hand side of equation~(\ref{eq:S}),
which amounts to an integration backwards in time. Integrating equation~(\ref{eq:S})
with the initial value $S_{ij}(0)=\hat{S}_{ij}$ until time $\tau$
yields $\tilde{S}_{ij}=S_{ij}(\tau)$.

In the third and last step, we set 
\[
P(t,x_{(1)},x_{(2)})=\sum_{i=1}^{r}X_{i}^{1}(t,x_{(1)})L_{i}(t,x_{(2)}),\quad\textrm{with}\quad L_{i}(t,x_{(2)})=\sum_{j=1}^{r}S_{ij}(t)X_{j}^{2}(t,x_{(2)})
\]
and use this representation for equation~(\ref{eq:LT_3}). Now $X_{i}^{1}$
remains constant and a similar calculation as for the two previous
steps yields the equation for the $L$ \emph{step}
\begin{equation}
\partial_{t}L_{i}(t,x_{(2)})=\sum_{\mu=1}^{M}\sum_{j=1}^{r}\left(c_{ij}^{2,\mu}(x_{(2)})L_{j}(t,x_{(2)}-\nu_{\mu,(2)})-d_{ij}^{2,\mu}(x_{(2)})L_{j}(t,x_{(2)})\right)\label{eq:L}
\end{equation}
with the time-independent coefficients 
\begin{equation}
\begin{aligned}c_{ij}^{2,\mu}(x_{(2)}) & =\langle X_{1,i}^{1}(x_{(1)}),a_{\mu}(x_{(1)}-\nu_{\mu,(1)},x_{(2)}-\nu_{\mu,(2)})X_{1,j}^{1}(x_{(1)}-\nu_{\mu,(1)})\rangle_{2},\\
d_{ij}^{2,\mu}(x_{(2)}) & =\langle X_{1,i}^{1}(x_{(1)}),a_{\mu}(x_{(1)},x_{(2)})X_{1,j}^{1}(x_{(1)})\rangle_{2}.
\end{aligned}
\label{eq:L_coeff}
\end{equation}
We integrate equation~(\ref{eq:L}) with the initial value 
\[
L_{i}(0,x_{(2)})=\sum_{j=1}^{r}\tilde{S}_{ij}X_{0,j}^{2}(x_{(2)})
\]
until time $\tau$ to obtain $L_{1,j}(x_{(2)})=L_{j}(\tau,x_{(2)})$.
Performing a QR decomposition 
\[
L_{1,i}(x_{(2)})=\sum_{j=1}^{r}S_{1,ij}X_{1,j}^{2}(x_{(2)})
\]
yields the orthonormal functions $X_{1,j}^{2}$ and the matrix $S_{1,ij}$,
which completes the first-order Lie-Trotter projector splitting algorithm.
The approximation to the solution at time $\tau$ is then given by
\[
P(\tau,x)\approx\sum_{i,j=1}^{r}X_{1,i}^{1}(x_{(1)})S_{1,ij}(t)X_{1,j}^{2}(x_{(2)}).
\]

We note that this approach can be easily extended to, e.g., the second-order
Strang splitting (see, e.g., \cite{Einkemmer_2018,Einkemmer_2023}).

\section{Algorithm and implementation\label{subsec:implementation-general-remarks}}

The CME can be regarded as an infinite system of ordinary differential
equations (ODEs) or as a discrete partial differential equation (PDE)
with spatial differences instead of derivatives \cite{Jahnke_2008}.
In order to turn the CME into a finite problem, we truncate the state
space to a finite domain which allows a numerical solution and still
captures enough of the information of the full (infinite) system.
If we define the truncated state space as $\Omega^{\zeta,\eta}=\{x\text{\ensuremath{\in\mathbb{N}_{0}^{N}}}:\zeta_{i}\le x_{i}\le\eta_{i}\ \mathrm{for}\ i=1,\dots,N\}$,
where $\zeta_{i}\in\mathbb{N}_{0}$ and $\eta_{i}\in\mathbb{N}_{0}$
and $\zeta_{i}<\eta_{i}$ ($i=1,\ldots,N)$, then the truncation error
can be estimated as follows:

We denote by $\mathscr{A}^{\zeta,\eta}$ the restriction of the linear
operator $\mathscr{A}$ (defined in equation \ref{eq:linear_operator})
to $\Omega^{\zeta,\eta}$ and by $P^{\zeta,\eta}(t)$ the solution
of the restricted CME $\partial_{t}P^{\zeta,\eta}(t)=\mathscr{A}^{\zeta,\eta}P^{\zeta,\eta}(t)$
with initial condition $P^{\zeta,\eta}(0)$, which is the initial
probability distribution restricted to the truncated state space.
Defining the total mass $m^{\zeta,\eta}=\sum_{x\in\Omega^{\zeta,\eta}}P^{\zeta,\eta}(t,x)$
and assuming that $m^{\zeta,\eta}\ge1-\epsilon$, \cite{Munsky_2006}
showed that
\[
P(t,x)-\epsilon\leq P^{\zeta,\eta}(t,x)\le P(t,x)\quad\text{for}\quad x\in\Omega^{\zeta,\eta}.
\]
These inequalities give an estimation of how close the truncated state
space solution approximates the true solution. The main issue of the
truncation is how to determine suitable $\zeta$ and $\eta$ for given
final time $t$ and tolerance $\epsilon>0$, such that $m^{\zeta,\eta}\ge1-\epsilon$.
Solving the reaction network deterministically with ODEs (which is
cheap) or biological insight into the system might give a good idea
on how to choose $\zeta$ and $\eta$ a priori. Alternatively, one
can implement a scheme with an adaptive truncated state space where
the error in mass is used as an indicator.

In our numerical implementation we work with truncated state spaces
$\Omega_{1}^{\zeta,\eta}=\{x_{(1)}\text{\ensuremath{\in\mathbb{N}_{0}^{m_{1}}}}:\zeta_{i}\le x_{i}\le\eta_{i}\ \mathrm{for}\ i=1,\dots,m_{1}\}$
and $\Omega_{2}^{\zeta,\eta}=\{x_{(2)}\text{\ensuremath{\in\mathbb{N}_{0}^{m_{2}}}}:\zeta_{i}\le x_{i}\le\eta_{i}\ \mathrm{for}\ i=m_{1}+1,\dots,N\}$
for the two partitions of the reaction network, where $\zeta_{i}$
and $\eta_{i}$ are fixed. We denote the number of degrees of freedom
by $n_{1}=(\zeta_{1}-\eta_{1}+1)\cdot\ldots\cdot(\zeta_{m_{1}}-\eta_{m_{1}}+1)$
and $n_{2}=(\zeta_{m_{1}+1}-\eta_{m_{1}+1}+1)\cdot\ldots\cdot(\zeta_{N}-\eta_{N}+1)$
for partition 1 and 2, respectively. The total number of degrees of
freedom is $n=n_{1}n_{2}$.

In the implementation we will store quantities depending on the population
number (such as $X_{i}^{1}(x_{(1)})$ and $X_{j}^{2}(x_{(2)})$) as
matrices. As $x_{(1)}$ and $x_{(2)}$ are vectors of size $m_{1}$
and $m_{2}$, respectively, we have to linearize the population number
dependency in order to store for example $X_{i}^{1}(x_{(1)})$ and
$X_{j}^{2}(x_{(2)})$ as matrices. We achieve this by introducing
bijective maps $\alpha:\,\Omega_{1}^{\zeta,\eta}\to\{1,\dots,n_{1}\}$
and $\beta:\,\Omega_{2}^{\zeta,\eta}\to\{1,\dots,n_{2}\}$ and thus
construct matrices $X^{1}=(\underline{X}_{1}^{1},\dots,\underline{X}_{r}^{1})\in\mathbb{R}^{n_{1}\times r}$
and $X^{2}=(\underline{X}_{1}^{2},\dots,\underline{X}_{r}^{2})\in\mathbb{R}^{n_{2}\times r}$,
whose columns are the low-rank factors evaluated on the truncated
state spaces $\Omega_{1}^{\zeta,\eta}$ and $\Omega_{2}^{\zeta,\eta}$.
Note that we indicate linearized quantities by underlining them, i.e.
$\underline{X}_{i}^{1}=(X_{i,\underline{\alpha}_{1}}^{1},\dots,X_{i,\underline{\alpha}_{n_{1}}}^{1})^{T}$,
where $\underline{\alpha}=(\alpha(x))_{x\in\Omega_{1}^{\zeta,\eta}}$.
The matrices $K\in\mathbb{R}^{n_{1}\times r}$and $L\in\mathbb{R}^{n_{2}\times r}$
are then computed by matrix multiplication, $K=X^{1}S$ and $L=X^{2}S^{T}$.

Using the substitution (\ref{eq:substitution}), the coefficients
(\ref{eq:K_coeff}) can be written as

\begin{equation}
\begin{aligned}C^{1,\mu}(x_{(1)}) & =\left[\mathcal{T}_{2,\mu}^{-1}[X_{0}^{2}]\right]^{T}\mathrm{diag}\left(\underline{a}_{\mu}(x_{(1)})\right)X_{0}^{2},\\
D^{1,\mu}(x_{(1)}) & =\left(X_{0}^{2}\right)^{T}\mathrm{diag}\left(\underline{a}_{\mu}(x_{(1)})\right)X_{0}^{2},
\end{aligned}
\label{eq:K_coeff_imp}
\end{equation}
with $X_{0}^{2}=X^{2}(t=0)\in\mathbb{R}^{n_{2}\times r}$ and $C^{1,\mu}(x_{(1)}),\,D^{1,\mu}(x_{(1)})\in\mathbb{R}^{r\times r}$.
The \emph{shift operator $\mathcal{T}_{2,\mu}=\mathcal{T}_{2,\mu}^{+1}$}
and the \emph{inverse shift operator }$\mathcal{T}_{2,\mu}^{-1}$
act element-wise and are defined as 
\[
\mathcal{T}_{2,\mu}^{\pm1}[X_{\underline{\beta}_{i}}^{2}]=\begin{cases}
0\qquad\mathrm{if} & \exists x_{j}=\left(\beta^{-1}(\underline{\beta}_{i})\right)_{j}:\,(x_{j}\pm\nu_{\mu,j}^{(2)}<\zeta_{j}^{(2)})\lor(x_{j}\pm\nu_{\mu,j}^{(2)}>\eta_{j}^{(2)}),\,j=1,\dots,m_{2}\\
X_{\underline{\beta}_{i}\pm\beta(\nu_{\mu}^{(2)})}^{2} & \mathrm{otherwise},
\end{cases}
\]
and $\zeta=(\zeta^{(1)},\zeta^{(2)})$ and $\eta=(\eta^{(1)},\eta^{(2)})$.
This definition approximates all terms $X^{2}(x_{(2)}-\nu_{\mu,(2)})$
which lie outside the truncated state space by 0, which assumes that
the probability function and the low rank factors have to decay sufficiently
fast within the truncated state space.

Writing $\underline{C}^{1,\mu}=(C_{\underline{\alpha}_{1}}^{1,\mu},\dots,C_{\underline{\alpha}_{n_{1}}}^{1,\mu})$,
the evolution equation of the $K$ step becomes 
\begin{equation}
\partial_{t}K=\left(\sum_{\mu=1}^{M}\mathcal{T}_{1,\mu}\left[\underline{K}\odot(\underline{C}^{1,\mu})^{T}\right]+\underline{K}\odot(\underline{D}^{1,\mu})^{T}\right),\label{eq:K_imp}
\end{equation}
with element-wise matrix-vector multiplication $\underline{K}\odot(\underline{D}^{1,\mu})^{T}=(K_{\underline{\alpha}_{1}}(D_{\underline{\alpha}_{1}}^{1,\mu})^{T},\dots,K_{\underline{\alpha}_{n_{1}}}(D_{\underline{\alpha}_{n_{1}}}^{1,\mu})^{T})$.
The shift operator $\mathcal{T}_{1,\mu}$ is defined in a similar
way as $\mathcal{T}_{2,\mu}$.

If we perform the integration over partition 2 in equations~(\ref{eq:S_coeff})
first, we can reuse the coefficients $C^{1,\mu}$ and $D^{1,\mu}$
for the calculation of the $S$ step coefficients
\begin{equation}
\begin{aligned}E_{ijkl} & =\sum_{\mu=1}^{M}\left(\mathcal{T}_{1,\mu}^{-1}\left[\underline{X}_{1,i}^{1}\right]\right)^{T}\mathrm{diag}\left(\underline{C}_{jl}^{1,\mu}\right)\underline{X}_{1,k}^{1},\\
F_{ijkl} & =\sum_{\mu=1}^{M}\left(\underline{X}_{1,i}^{1}\right)^{T}\mathrm{diag}\left(\underline{D}_{jl}^{1,\mu}\right)\underline{X}_{1,k}^{1},
\end{aligned}
\label{eq:S_coeff_imp}
\end{equation}
with $E^{\mu},\,F^{\mu}\in\mathbb{R}^{r\times r\times r\times r}$.
With these coefficients we can write the evolution equation of the
$S$ step as
\begin{equation}
\partial_{t}S_{ij}=-\sum_{k,l=1}^{r}S_{kl}\left(E_{ijkl}-F_{ijkl}\right).\label{eq:S_imp}
\end{equation}
The coefficients $C^{2,\mu}$, $D^{2,\mu}$ are calculated via 
\begin{equation}
\begin{aligned}C^{2,\mu}(x_{(2)}) & =\left[\mathcal{T}_{1,\mu}^{-1}[X_{1}^{1}]\right]^{T}\mathrm{diag}\left(\underline{a}_{\mu}(x_{(2)})\right)X_{1}^{1},\\
D^{2,\mu}(x_{(2)}) & =\left(X_{1}^{1}\right)^{T}\mathrm{diag}\left(\underline{a}_{\mu}(x_{(2)})\right)X_{1}^{1},
\end{aligned}
\label{eq:L_coeff_imp}
\end{equation}
with $X_{1}^{1}=X^{1}(t=\tau)\in\mathbb{R}^{n_{1}\times r}$ and $C^{2,\mu}(x_{(2)}),\,D^{2,\mu}(x_{(2)})\in\mathbb{R}^{r\times r}$. 

Showing a similar structure as the corresponding equation for the
$K$ step, the evolution equation for the $L$ step reads as
\begin{equation}
\partial_{t}L=\sum_{\mu=1}^{M}\left(\mathcal{T}_{2,\mu}\left[\underline{L}\odot(\underline{C}^{2,\mu})^{T}\right]+\underline{L}\odot(\underline{D}^{2,\mu})^{T}\right).\label{eq:L_imp}
\end{equation}
Note that for the second term on the right-hand side of equation~(\ref{eq:L_imp})
we could perform the summation over all reactions $R_{\mu}$ before
multiplying $\underline{D}^{2,\mu}$ with $\underline{L}$. Moreover,
the calculation of the evolution equation~(\ref{eq:L_imp}) could
be simplified by introducing a reaction-independent $\underline{D}^{2}=\sum_{\mu=1}^{M}\underline{D}^{2,\mu}$.
However, as we will see in section~\ref{subsec:reagents_dependence},
it is computationally more efficient to perform the summation over
all reactions after multiplying $\underline{L}$ with the reaction-dependent
$\underline{D}^{2,\mu}$. The same holds for the coefficient $\underline{D}^{1,\mu}$
and the second term on the right-hand side of the evolution equation
(\ref{eq:K_imp}) for the $K$ step. For the first terms on the right-hand
side of equations~(\ref{eq:K_imp}) and (\ref{eq:L_imp}) we always
have to keep the reaction-dependence of $\underline{C}^{1,\mu}$ and
$\underline{C}^{2,\mu}$, since the shift operator is also reaction-dependent,
therefore the computational effort would scale the same even when
using reaction-independent $D$ coefficients. 

Finally, we want to give a remark on the computational effort for
the evolution equation and the calculation of the coefficients. Without
making any further simplifications, the computational effort for calculating
the $C$ and $D$ coefficients is $\mathcal{O}(Mr^{2}n)$, where $M$
was the total number of reactions channels. When we reuse $\underline{C}^{1,\mu}$
and $\underline{D}^{1,\mu}$ for the calculation of the $E$ and $F$
coefficients, the complexity for computing the $E$ and $F$ coefficients
is $\mathcal{O}(M\,r^{4}n_{1})$. The right-hand side of the evolution
equation for the $K$ step requires an integration over the population
numbers in partition 1, therefore the computational cost is $\mathcal{O}(Mr^{2}n_{1})$.
Similarly, the computation of the $L$ step scales with $\mathcal{O}(Mr^{2}n_{2})$,
whereas for the $S$ step we have complexity $\mathcal{O}(r^{4})$.
Note that in particular the computation of the $C$ and $D$ coefficients
is very expensive, since it scales with the total number of degrees
of freedom $n$. Thus it is imperative to reduce this computational
burden, which is the topic of the next section.

\subsection{Efficient computation of the coefficients\label{subsec:efficient-computation}}

In the previous section we have seen that computing the coefficients
without making any further assumptions is computationally expensive.
The main goal here is to describe ways how to avoid the scaling of
the computational effort with the total number of degrees of freedom
$n$. We essentially discuss two possibilities to circumvent this
scaling behaviour: First, most reactions only depend on a small subset
of all species. When we denote for a given reaction $R_{\mu}$ the
number of participating species (we call them \emph{reagents}) by
$\tilde{N}_{\mu}$, than this assumption can be expressed as $\tilde{N}_{\mu}\ll N$.
Second, in many reaction networks the propensity functions exhibit
a so-called factorization property, and exploiting this property again
reduces the computational burden. Note that our present implementation
does not exploit the factorization property since in all our examples
$\tilde{N}_{\mu}\ll N$. However, what the discussion in this section
shows is that even in the rare instances where this is not the case,
the factorization property which is common to most reactions gives
a way forward to efficiently implementing the dynamical low-rank approach. 

\subsubsection{Dependence of the propensity on reagents\label{subsec:reagents_dependence}}

Since in most reactions only a subset of all species is actually participating,
the propensity $a_{\mu}(x)$ for such a reaction $R_{\mu}$ only depends
on the population number of the $\tilde{N}_{\mu}$ reagents, so $a_{\text{\ensuremath{\mu}}}(x)=a_{\mu}(\tilde{x}_{\mu})$,
where $\tilde{x}_{\mu}\in\mathbb{N}_{0}^{\tilde{N}_{\mu}}$ and $\tilde{N}_{\mu}\le N$.
In many cases the propensities only depend on the population number
of two or three species, therefore $\tilde{N}_{\mu}\ll N$. The computational
effort of our algorithm can be reduced substantially by calculating
coefficients $C$ and $D$ only for the possible values of $\tilde{x}_{\mu}$
that are actually needed.

For a given reaction $\mu$ we first determine in the implementation
the reagents and precompute all possible values of the propensity
function $a_{\mu}(\tilde{x}_{\mu})$, since those values do not change
over time. The $C$ and $D$ coefficients have to be calculated only
for $\tilde{n}_{1}^{\mu}$ or $\tilde{n}_{2}^{\mu}$ population number
values, but $K$ and $L$ still depend on the population numbers $x_{(1)}$
and $x_{(2)}$, respectively. Therefore we have to introduce a (reaction-dependent)
mapping between $\tilde{x}_{\mu}$ and $x_{(1)}$ in order to perform
for example the multiplication $\underline{K}\odot(\underline{D}^{1,\mu})^{T}$
on the right-hand side in equation~(\ref{eq:K_imp}). Applying this
map effectively introduces a reaction-dependency on the overall multiplication
term. This is the reason why we cannot introduce the reaction independent
$\underline{D}^{2}=\sum_{\mu=1}^{M}\underline{D}^{2,\mu}$ as discussed
previously.

Note that the complexity for the integration over $x_{1}$ in equation~(\ref{eq:K_coeff_imp})
still scales with $\mathcal{O}(n_{1})$ and for equation~(\ref{eq:L_coeff_imp})
the integration over $x_{2}$ scales with $\mathcal{O}(n_{2})$, but
the coefficients have to be calculated only for the $\tilde{n}_{1}^{\mu}$
or $\tilde{n}_{2}^{\mu}$ population number values. Therefore the
complexity for calculating the $\underline{C}^{1,\mu}$ and $\underline{D}^{1,\mu}$
coefficients is reduced to $\mathcal{O}(\sum_{\mu=1}^{M}\tilde{n}_{1}^{\mu}n_{2}r^{2})$,
and for the $\underline{C}^{2,\mu}$ and $\underline{D}^{2,\mu}$
coefficients to $\mathcal{O}(\sum_{\mu=1}^{M}\tilde{n}_{2}^{\mu}n_{1}r^{2})$.
Thus, these computations do no longer scale with $n$ (assuming that
$\tilde{N}_{\mu}\ll N$). 

\subsubsection{Factorization property of the propensity function}

The equations for the coefficients (\ref{eq:K_coeff}), (\ref{eq:S_coeff})
and (\ref{eq:L_coeff}) can be simplified if the propensity function
can be written as 
\begin{equation}
a_{\mu}(x_{(1)},x_{(2)})=a_{\mu,(1)}(x_{(1)})\,a_{\mu,(2)}(x_{(2)}).\qquad\textrm{(factorization property)}\label{eq:factor}
\end{equation}
This property is valid for elementary reaction types and reactions
of the Michaelis-Menten form\textcolor{black}{{} and thus is ubiquitous
in most biological systems.} The factorization property enables us
to rewrite for example the coefficient $c_{ij}^{1,\mu}$ in equation~(\ref{eq:K_coeff})
as 
\[
c_{ij}^{1,\mu}(x_{(1)})=a_{\mu,(1)}(x_{(1)}-\nu_{\mu,(1)})\langle X_{0,i}^{2}(x_{(2)}),a_{\mu,(2)}(x_{(2)}-\nu_{2}^{\mu})X_{0,j}^{2}(x_{2}-\nu_{2}^{\mu})\rangle_{2},
\]
which scales with $\mathcal{O}(M\,r^{2}(n_{1}+n_{2})$ compared to
$\mathcal{O}(M\,r^{2}n)$ (even when disregarding the considerations
about the dependence of the propensity on reagents in the previous
section). Moreover, the two inner products of the coefficient $e_{ijkl}$
in equation~(\ref{eq:S_coeff}) can be calculated independently,
\[
\begin{split}e_{ijkl} & =\sum_{\mu=1}^{M}\langle X_{1,i}^{1}(x_{(1)})a_{\mu,(1)}(x_{(1)}-\nu_{\mu,(1)})X_{1,k}^{1}(x_{(1)}-\nu_{\mu,(1)})\rangle_{1}\\
 & \quad\times\langle X_{0,j}^{2}(x_{(2)})a_{\mu,(2)}(x_{(2)}-\nu_{\mu,(2)})X_{0,l}^{1}(x_{(2)}-\nu_{\mu,(2)})\rangle_{2},
\end{split}
\]
which has computational costs of $\mathcal{O}(Mr^{4}(n_{1}+n_{2}))$
instead of $\mathcal{O}(Mr^{4}n)$.

\subsection{First- and second-order projector splitting integrator}

The first-order integrator is obtained by using Lie--Trotter splitting
as explained in section~\ref{sec:DLR-approximation}. The evolution
equations (\ref{eq:K_imp}) and (\ref{eq:L_imp}) for $K$ and $L$
steps as well as (\ref{eq:S_imp}) for the $S$ step are solved with
an explicit Euler method. Due to different reaction time scales stemming
from both small and large propensity values, the CME becomes stiff
for many systems. In order to remain in the stable region for large
time step size $\tau$, we perform $k$ explicit Euler steps with
a time step size of $\tau/k$ while keeping the coefficients constant.
Note that although the computational cost for evaluating the right-hand
side of the evolution equation shows the same scaling, the constant
is smaller compared to the calculation of the coefficients, which
however only needs to be done once or twice (for $C^{1,\mu}$ and
$D^{1,\mu}$ in case of the second-order integrator) in each time
step. We will explore the use of implicit integrators in future work.

The low-rank factors $X^{1}$ and $X^{2}$ and the coupling coefficients
$S$ are obtained from $K$ and $L$ matrices by performing a QR decomposition.
In order to perform the QR decomposition and the linear algebra operations
required for an efficient calculation of the coefficients we made
use of the dynamical low-rank framework \texttt{Ensign}~\cite{Cassini_2021}.

A detailed description of the first order Lie-Trotter projector splitting
scheme is shown in algorithm \ref{alg:first-order}.
\begin{algorithm}[H]
\caption{\label{alg:first-order}First-order Lie--Trotter projector splitting
integrator for the kinetic CME.}

\textbf{Input:} $X_{0}^{1},$ $S_{0}$, $X_{0}^{2}$

\textbf{Output:} $X_{1}^{1}$, $S_{3}$, $X_{1}^{2}$

\begin{algorithmic}[1]
\State Calculate $C^{1,\mu}(x_{(1)})$ and $D^{1,\mu}(x_{(1)})$ with $X_0^2$ using equation~(\ref{eq:K_coeff_imp})
\State Integrate $K$ from $0$ to $\tau$ with initial value $K(0) = X_0^1 S_0$ using equation~(\ref{eq:K_imp})
\State Decompose $K(\tau) = X_1^1 S_1$ via a QR factorization
\State Calculate $E^{\mu}$ and $F^{\mu}$ with $X_1^1$, $X_0^2$, $C^{1,\mu}(x_{(1)})$ and $D^{1,\mu}(x_{(1)})$ using equation~(\ref{eq:S_coeff_imp})
\State Integrate $S$ from $0$ to $\tau$ with initial value $S(0) = S_1$ using equation~(\ref{eq:S_imp}) and set $S_2 = S(\tau)$
\State Calculate $C^{2,\mu}(x_{(2)})$ and $D^{2,\mu}(x_{(2)})$ with $X_1^1$ using equation~(\ref{eq:L_coeff_imp})
\State Integrate $L$ from $0$ to $\tau$ with initial value $L(0) = X_1^2 (S_2)^T$ using equation~(\ref{eq:L_imp})
\State Decompose $L(\tau) = X_1^2 (S_3)^T$ via a QR factorization
\end{algorithmic}
\end{algorithm}

Our numerical scheme can be generalized to a second-order method by
employing Strang splitting in the context of equation~(\ref{eq:projector-equation}).
This is shown in detail in algorithm \ref{alg:second-order}. Note
that two of the steps are repeated while one step is only performed
once (due to the symmetry of the splitting). Ideally, the step that
has to be done only once is chosen to coincide with the step that
incurs the largest computational effort (either the $K$ or the $L$
step).

\begin{algorithm}[H]
\caption{\label{alg:second-order}Second-order Strang projector splitting integrator
for the kinetic CME.}

\textbf{Input:} $X_{0}^{1},$ $S_{0}$, $X_{0}^{2}$

\textbf{Output:} $X_{2}^{1}$, $S_{5}$, $X_{1}^{2}$

\begin{algorithmic}[1]
\State Calculate $C^{1,\mu}(x_{(1)})$ and $D^{1,\mu}(x_{(1)})$ with $X_0^2$ using equation~(\ref{eq:K_coeff_imp})
\State Integrate $K$ from $0$ to $\tau/2$ with initial value $K(0) = X_0^1 S_0$ using equation~(\ref{eq:K_imp})
\State Decompose $K(\tau/2) = X_1^1 S_1$ via a QR factorization
\State Calculate $E^{\mu}$ and $F^{\mu}$ with $X_1^1$, $X_0^2$, $C^{1,\mu}(x_{(1)})$ and $D^{1,\mu}(x_{(1)})$ using equation~(\ref{eq:S_coeff_imp})
\State Integrate $S$ from $0$ to $\tau/2$ with initial value $S(0) = S_1$ and set $S_2 = S(\tau/2)$ using equation~(\ref{eq:S_imp})
\State Calculate $C^{2,\mu}(x_{(2)})$ and $D^{2,\mu}(x_{(2)})$ with $X_1^1$ using equation~(\ref{eq:L_coeff_imp})
\State Integrate $L$ from $0$ to $\tau$ with initial value $L(0) = X_1^2 (S_2)^T$ using equation~(\ref{eq:L_imp})
\State Decompose $L(\tau) = X_1^2 (S_3)^T$ via a QR factorization
\State Recalculate $C^{1,\mu}(x_{(1)})$ and $D^{1,\mu}(x_{(1)})$ with $X_1^2$  using equation~(\ref{eq:K_coeff_imp})
\State Recalculate $E^{\mu}$ and $F^{\mu}$ with $X_1^1$, $X_1^2$ and new values for $C^{1,\mu}(x_{(1)})$ and $D^{1,\mu}(x_{(1)})$ using equation~(\ref{eq:S_coeff_imp})
\State Integrate $S$ from $\tau/2$ to $\tau$ with initial value $S(\tau/2) = S_3$ and set $S_4 = S(\tau)$ using equation~(\ref{eq:S_imp})
\State Integrate $K$ from $\tau/2$ to $\tau$ with initial value $K(\tau/2) = X_1^1 S_4$ using equation~(\ref{eq:K_imp})
\State Decompose $K(\tau) = X_2^1 S_5$ via a QR factorization
\end{algorithmic}
\end{algorithm}

\section{Numerical experiments\label{sec:numerical-experiments}}

We tested our implementation with three models from the field of biochemistry.
The smallest model, the genetic toggle switch, was primarily chosen
for code validation and to investigate the approximation accuracy
(as a reference solution without the low-rank approximation can be
computed easily). For the the two larger models, the bacteriophage-$\lambda$
(``lambda phage'') and the BAX pore assembly, we compare the DLR
approximation with the dominating numerical method for solving the
CME, the stochastic simulation algorithm (SSA) (see, e.g., \cite{Gillespie_1976}).

\subsection{Toggle switch}

The genetic toggle switch, as first described in \cite{Gardner_2000},
has a function analogous to a flip-flop in electronics. It consists
of two mutually repressing proteins $S_{1}$ and $S_{2}$, which leads
to two stable steady-states. We studied the reaction system shown
in table~\ref{tab:ts}, which was also considered in \cite{Jahnke_2008}.
\begin{table}[H]
\begin{centering}
\begin{tabular}{c|c|c}
\hline 
No. & Reaction & Propensity function\tabularnewline
\hline 
1 & $S_{1}\longrightarrow\star$ & $c\cdot x_{1}$\tabularnewline
2 & $S_{2}\longrightarrow\star$ & $c\cdot x_{2}$\tabularnewline
3 & $\star\longrightarrow S_{1}$ & $b/(b+x_{2})$\tabularnewline
4 & $\star\longrightarrow S_{2}$ & $b/(b+x_{1})$\tabularnewline
\hline 
\end{tabular}
\par\end{centering}
\caption{Reactions and propensity functions of the toggle switch systems. The
two parameters are chosen as $b=0.4$ and $c=0.05$.\label{tab:ts}}
\end{table}
The first two reactions describe the decay of proteins $S_{1}$ and
$S_{2}$, respectively. If the population number of $S_{2}$ is large,
then the propensity of reaction 3 becomes small and transcription
of new copies of $S_{1}$ is inhibited. Similarly, the production
of $S_{2}$ by reaction 4 is inhibited by $S_{1}$.

As initial value we consider the Gaussian distribution
\begin{align*}
P(0,x) & =\gamma\cdot\exp\left(-\frac{1}{2}(x-\mu)^{T}C^{-1}(x-\mu)\right),\\
C & =\frac{1}{2}\begin{pmatrix}75 & -15\\
-15 & 75
\end{pmatrix},
\end{align*}
with $\mu=(30,\,5)$ and $\gamma$ was determined by the condition
$\sum_{x\in\Omega^{\zeta,\eta}}P(0,x)=1$.

We solved the CME on the time interval $[0,500]$ with truncation
indices $\eta=(0,0)$ and $\zeta=(50,50)$ and with the trivial partitions
$\mathcal{P}_{1}=\{S_{1}\}$ and $\mathcal{P}_{2}=\{S_{2}\}$. Using
rank $r=5$, the total number of degrees of freedom is reduced from
$51^{2}=2601$ to $2\cdot51\cdot5+5^{2}=535$, which is $20.6\%$
of the full system size. Due to the relatively small size of the truncated
state space an ``exact'' reference solution of the full system on
the truncated state space could be obtained via a Python implementation
that uses the RK45 \texttt{scipy.solve\_ivp} routine to directly solve
equation (\ref{eq:CME}).

\begin{figure}[H]
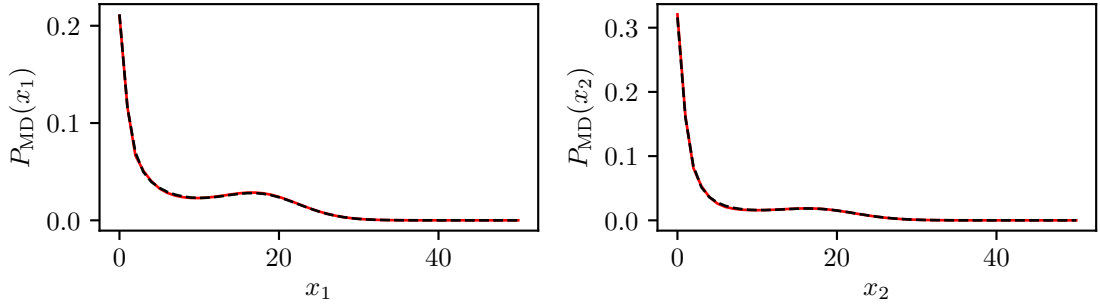

\centering
\input{plots/ts_fig1_r4.pgf}\\
\input{plots/ts_fig1_r5.pgf}

\caption{DLR approximation (red, solid line) and exact reference solution (black,
dashed line) of the one-dimensional marginal distributions $P_{\mathrm{MD}}(x_{1})$
and $P_{\mathrm{MD}}(x_{2})$ for ranks $r=4$ and $5$ of the toggle
switch example at $t=500$. The reference solution was obtained by
solving the CME directly with the \texttt{scipy.solve\_ivp} routine.
For the DLR approximation the second-order integrator with time step
size $\tau=0.02$ and $10$ substeps were used.\label{fig:ts1}}
\end{figure}

\begin{figure}[H]
\centering
\input{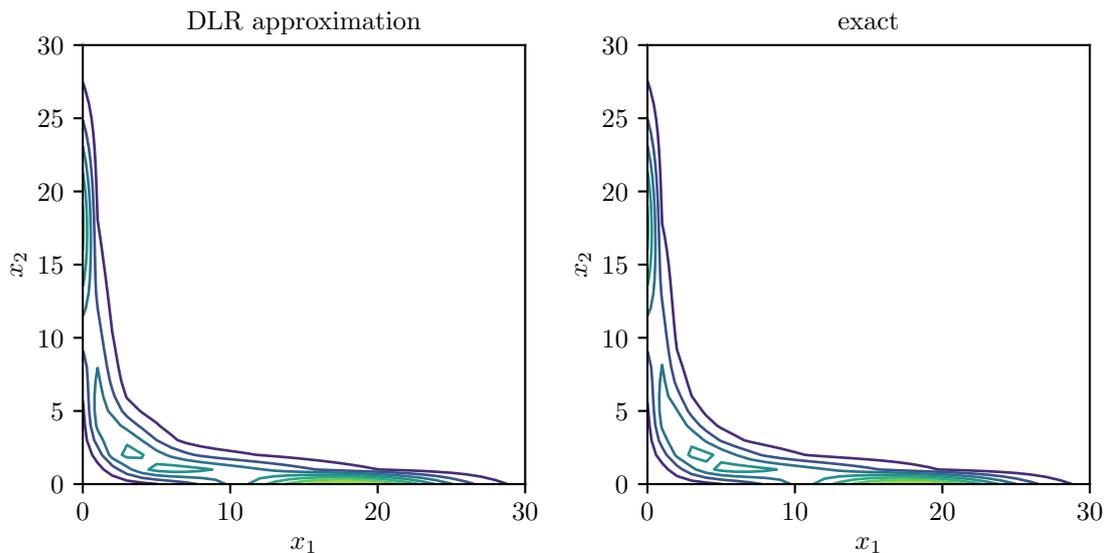}

\caption{DLR approximation (left) and exact reference solution (right) of the
full probability distribution $P(x_{1},x_{2})$ for the toggle switch
example at $t=500$. Note that $P$ is only defined at the discrete
grid points $x\in\mathbb{N}_{0}^{2}$; the contour plots are based
on interpolation and are shown here for the sake of clarity. The reference
solution was obtained by solving the CME directly with the \texttt{scipy.solve\_ivp}
routine. For the DLR approximation rank $r=5$ and the second-order
integrator with time step size $\tau=0.02$ and $10$ substeps were
used.\label{fig:ts2}}
\end{figure}

\begin{figure}[H]
\centering
\input{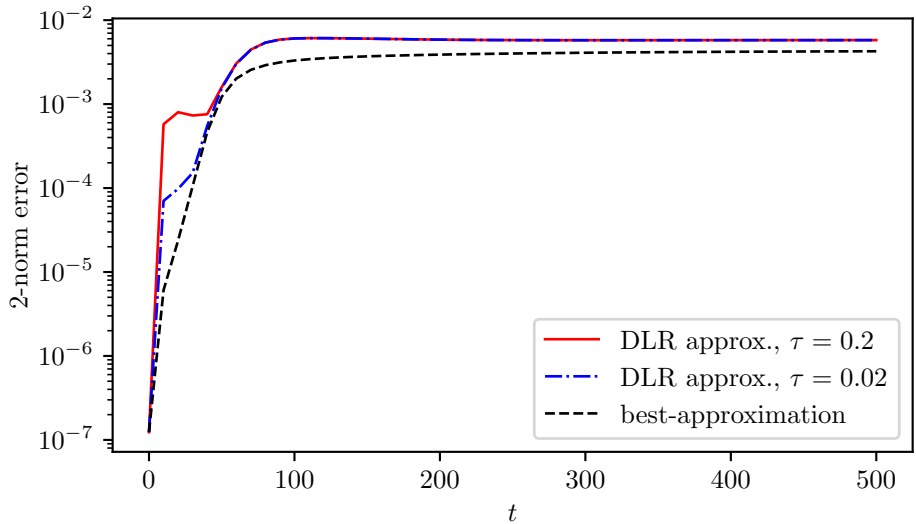}

\caption{Comparison of the 2-norm error for the DLR approximation using time
step size $\tau=0.2$ (red, solid line) and $\tau=0.02$ (blue, dash-dotted
line) with the best-approximation (black, dashed line) as a function
of time $t$ for the toggle switch system. The error was calculated
by comparing each approximation with the exact reference solution
of the toggle switch. The best-approximation was obtained by truncating
all but the first $r=5$ singular values of a SVD of the reference
solution. For the DLR approximation rank $r=5$ and the second-order
integrator with $10$ substeps were used.\label{fig:ts3}}
\end{figure}

Figure~\ref{fig:ts1} shows the DLR approximation (using the second-order
integrator with time step size $\tau=0.02$ and $10$ substeps) with
ranks $r=4$ and $5$ and the reference solution of the one-dimensional
marginal distributions $P_{\mathrm{MD}}(x_{1})$, $P_{\mathrm{MD}}(x_{2})$
at time $t=500$. It can be clearly seen that $r=4$ is not sufficient
to capture the full behavior of the system, but for $r=5$ we obtain
very good results. Figure~\ref{fig:ts2} depicts the full probability
distribution $P(x_{1},x_{2})$ at time $t=500$. This figure again
demonstrates that the results of the DLR approximation for $r=5$
are in very good agreement with the exact solution of the truncated
CME. Using the second-order integrator with time step size $\text{\ensuremath{\tau=0.02}}$
and $10$ substeps, the total run time for the simulation with rank
$r=5$ was approximately $1$ minute and $16$ seconds on a MacBook
Pro with a $2$ GHz Intel Core i5 Skylake (6360U) processor. The results
of the DLR approximation were computed with one thread.
Figure~\ref{fig:ts3} shows the $2$-norm error of the best-approximation
and of the DLR approximation for time step sizes $\tau=0.2$ and $0.02$,
using the second-order integrator with $10$ substeps. The best-approximation
was obtained by truncating all but the first $r=5$ singular values
of a singular value decomposition (SVD) of the reference solution,
for the DLR approximation we again used $r=5$. It can be seen that
using a smaller time step size helps particularly in the first few
steps of the simulation. After approximately $t=50$ the errors for
the different time step sizes are almost identical, indicating that
the overall error is dominated by the low-rank approximation. The
error of the dynamical low-rank algorithm proposed is only slightly
larger than the theoretical best approximation with the same rank.

\subsection{Lambda phage}

As a second example, the DLR approximation was applied to the model
for the life cycle of the lambda phage as described in \cite{Hegland_2007}.
Table~\ref{tab:lp} lists the ten reactions and five species of this
system.
\begin{table}[H]
\centering{}%
\begin{tabular}{c|c|c}
\hline 
No. & Reaction & Propensity function\tabularnewline
\hline 
1 & $\star\longrightarrow S_{1}$ & $a_{1}b_{1}/(b_{1}+x_{2})$\tabularnewline
2 & $\star\longrightarrow S_{2}$ & $(a_{2}+x_{5})b_{2}/(b_{2}+x_{1})$\tabularnewline
3 & $\star\longrightarrow S_{3}$ & $a_{3}b_{3}x_{2}/(b_{3}x_{2}+1)$\tabularnewline
4 & $\star\longrightarrow S_{4}$ & $a_{4}b_{4}x_{3}/(b_{4}x_{3}+1)$\tabularnewline
5 & $\star\longrightarrow S_{5}$ & $a_{5}b_{5}x_{3}/(b_{5}x_{3}+1)$\tabularnewline
6 & $S_{1}\longrightarrow\star$ & $c_{1}\cdot x_{1}$\tabularnewline
7 & $S_{2}\longrightarrow\star$ & $c_{2}\cdot x_{2}$\tabularnewline
8 & $S_{3}\longrightarrow\star$ & $c_{3}\cdot x_{3}$\tabularnewline
9 & $S_{4}\longrightarrow\star$ & $c_{4}\cdot x_{4}$\tabularnewline
10 & $S_{5}\longrightarrow\star$ & $c_{5}\cdot x_{5}$\tabularnewline
\hline 
\end{tabular}\qquad{}%
\begin{tabular}{c|c|c|c|c|c}
\hline 
 & $i=1$ & $i=2$ & $i=3$ & $i=4$ & $i=5$\tabularnewline
\hline 
$a_{i}$ & $0.5$ & $1$ & $0.15$ & $0.3$ & $0.3$\tabularnewline
$b_{i}$ & $0.12$ & $0.6$ & $1$ & $1$ & $1$\tabularnewline
$c_{i}$ & $0.0025$ & $0.0007$ & $0.0231$ & $0.01$ & $0.01$\tabularnewline
\hline 
\end{tabular}\caption{Reactions, propensity functions and parameters of the lambda phage
system.\label{tab:lp}}
\end{table}
The life cycle of the lambda phage represents a naturally occurring
toggle switch. The lambda phage infects \emph{E. coli}, and depending
on the environment, either stays dormant in the bacterial host (\emph{lysogenic
phase}) or multiplies, reassembles itself and breaks out of the host
(\emph{lytic phase}). If enough $S_{5}$ is present in the environment,
$S_{2}$ is produced and the system is in the lysogenic phase. Abundance
of $S_{2}$ in turn inhibits the formation of $S_{1}$ via reaction
1. If the amount of $S_{5}$ in the environment is scarce, the production
of $S_{1}$ causes the system to enter the lytic phase and the transcription
of new copies of $S_{2}$ via reaction 2 is inhibited.

As an initial value the multinomial distribution with parameters $n=3$
and $p=(0.05,\dots,0.05)$ has been chosen:
\[
P(0,x)=\begin{cases}
\frac{3!}{x_{1}!\cdots x_{5}!(3-|x|)!}0.05^{|x|}(1-5\cdot0.05)^{3-|x|} & \text{if}\quad|x|\le3,\\
0 & \text{else,}
\end{cases}
\]
where $|x|=x_{1}+\dots+x_{5}$. We solved the CME on the time interval
$[0,10]$ with truncation indices $\eta=(0,0,0,0,0)$ and $\zeta=(15,40,10,10,10)$.
The reaction network was partitioned into $\mathcal{P}_{1}=\{S_{1},S_{2}\}$
and $\mathcal{P}_{2}=\{S_{3},S_{4},S_{5}\}$, therefore the two partitions
have a comparable number of degrees of freedom, $n_{1}=16\cdot41=656$
and $n_{2}=11^{3}=1331$. Using rank $r=9$, the total number of degrees
of freedom used in the DLR approximation is reduced from $n_{1}\cdot n_{2}=873\,136$
to $(n_{1}+n_{2})\cdot r+r^{2}=17\,964$, which is $2.1\%$ of the
full system size.

An ``exact'' reference solution was obtained again by solving the
full CME on the truncated state space with \texttt{scipy.solve\_ivp}.
Due to the relatively large system size the computation of the full
solution is very costly and therefore substantially slower than the
DLR approximation. Moreover, we compare the DLR approximation with
results obtained with SSA. These latter results were computed in the
systems biology framework \emph{PySB} \cite{Lopez_2013}, which uses
the SSA implementation \emph{StochKit2} \cite{Sanft_2011}. Table~\ref{tab:lp_runtime}
gives an overview of the run times for the exact reference solution,
the DLR approximation and SSA.
\begin{table}[H]
\begin{centering}
\begin{tabular}{c|r}
\hline 
 & run time {[}s{]}\tabularnewline
\hline 
DLR approx. ($r=4$) & $52$\tabularnewline
DLR approx. ($r=9$) & $191$\tabularnewline
SSA ($10\,000$ runs) & $855$\tabularnewline
SSA ($100\,000$ runs) & $905$\tabularnewline
SSA ($1\,000\,000$ runs) & $2319$\tabularnewline
exact & $1164$\tabularnewline
\hline 
\end{tabular}
\par\end{centering}
\caption{Overview of the approximate run times in seconds for the DLR approximation,
SSA, and the exact reference solution for the lambda phage system.
The DLR approximation was computed with the second-order integrator
using time step size $\tau=0.01$ and $10$ substeps. All computations
were performed on a MacBook Pro with a $2$ GHz Intel Core i5 Skylake
(6360U) processor. The results of the DLR approximation were computed
with one thread.\label{tab:lp_runtime}}
\end{table}

\begin{figure}[H]
\centering
\input{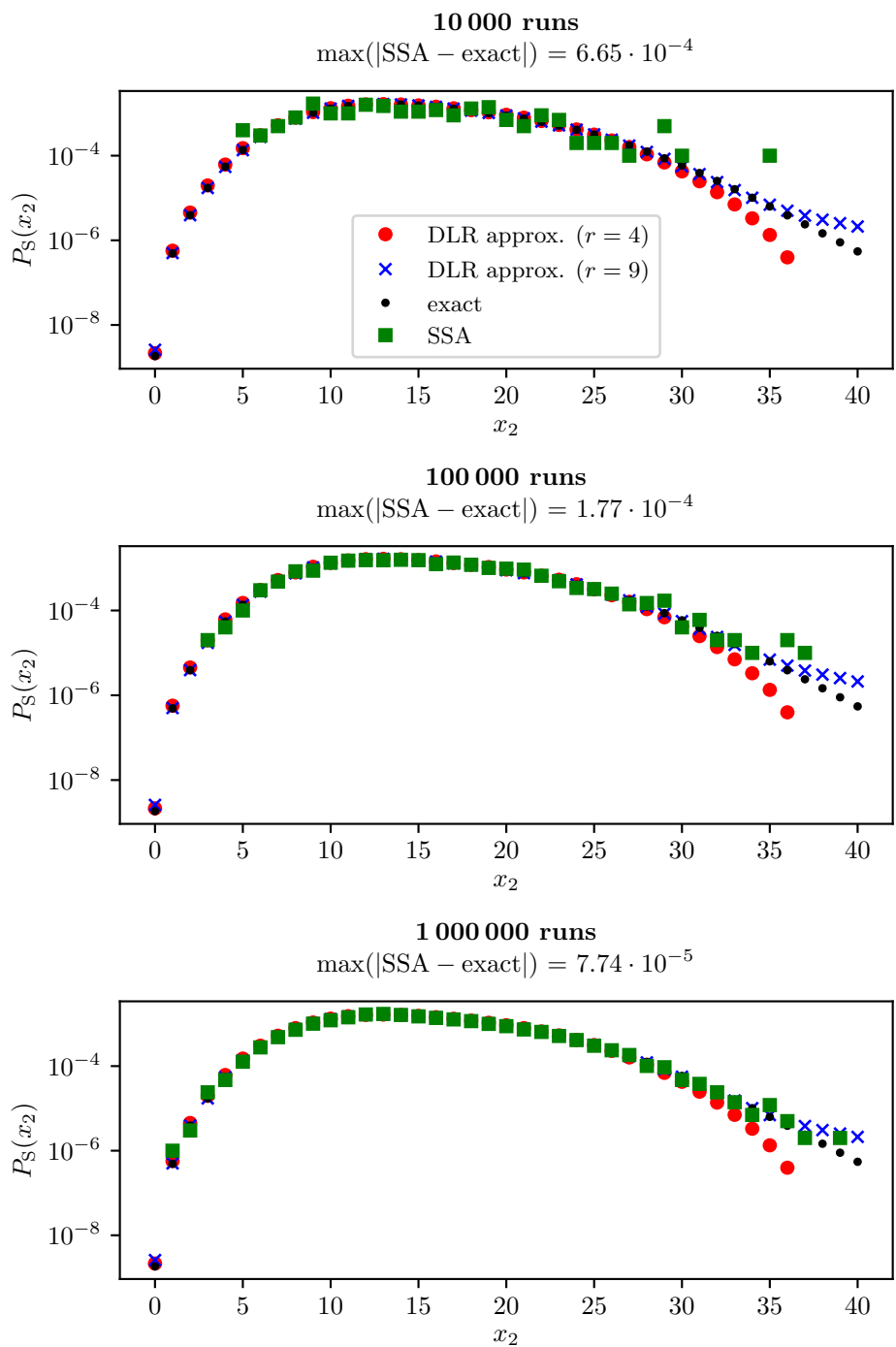}

\caption{Partially evaluated probability distribution $P_{\mathrm{S}}(x_{2})=P(x_{1}=0,x_{2},x_{3}=1,x_{4}=1,x_{5}=1)$
of the lambda phage system at $t=10$. The results were obtained with
our implementation of the DLR approximation using rank $r=4$ (red
dots), rank $r=9$ (blue crosses) and with SSA (green squares) using
$10\,000$, $100\,000$ and $1\,000\,000$ samples. For comparison
also the exact reference solution is shown (small black dots), which
was obtained by solving the CME on the truncated state space directly
with the \texttt{scipy.solve\_ivp} routine. The DLR approximation
was computed with the second-order integrator using time step size
$\tau=0.01$ and $10$ substeps. The maximum error for the DLR approximation
was $\max(|\mathrm{DLR}-\mathrm{exact}|)=4.21\cdot10^{-5}$ for $r=4$
and $1.03\cdot10^{-5}$ for $r=9$.\label{fig:lp1}}
\end{figure}

SSA is a Monte Carlo approach, therefore the results of this method
are polluted with noise which scales only as the inverse square root
of the total number of independent runs or samples. Figure~\ref{fig:lp1}
shows the partially evaluated probability distribution $P_{\mathrm{S}}(x_{2})=P(x_{1}=0,x_{2},x_{3}=1,x_{4}=1,x_{5}=1)$
at time $t=10$ computed with the DLR approximation using rank $r=4$
and $r=9$, as well as with SSA using $10\,000$, $100\,000$ and
$1\,000\,000$ samples. Furthermore, the exact reference solution
of the CME on the truncated state space is shown for comparison. The
DLR solution for $r=4$ is in very good agreement with the reference
solution, only for high population numbers $x_{2}>30$, which have
a relatively low probability, a discrepancy becomes visible. For population
numbers $x_{2}>36$ the results are still close to zero, but become
negative and therefore are not shown in this semi-logarithmic plot.
When increasing the rank to $r=9$, the results for $x_{2}>36$ remain
positive and only show a small deviation from the reference solution.
The results for SSA with $10\,000$ runs exhibit a lot of noise and
for several population numbers the probability is zero, because the
corresponding state was not sampled at all during the simulation.
When performing the simulation with more runs, this stochastic noise
decreases, but even for $1\,000\,000$ runs SSA has still problems
to resolve the small values at the tails of the probability distribution
(for example at $x_{2}=0$). This, in particular, shows that the dynamical
low-rank approximation has a significant advantage if one is interested
in resolving states with low probability. The reason for this is that
no numerical noise is introduced by the low-rank approach. We also
calculated the maximal error between the exact solution and the SSA
result for $P_{\mathrm{S}}(x_{2})$. For $1\,000\,000$ runs the maximal
error was $7.74\cdot10^{-5}$ and therefore higher than for the DLR
approximation for both $r=4$ ($4.21\cdot10^{-5}$) and $r=9$ ($1.03\cdot10^{-5}$).
The computational cost of the DLR approximation, despite the lower
error, is lower by a factor of approximately $45$ ($r=4$) and $12$
($r=9$) compared to SSA with $1\,000\,000$ samples. 

\begin{figure}[H]
\centering
\input{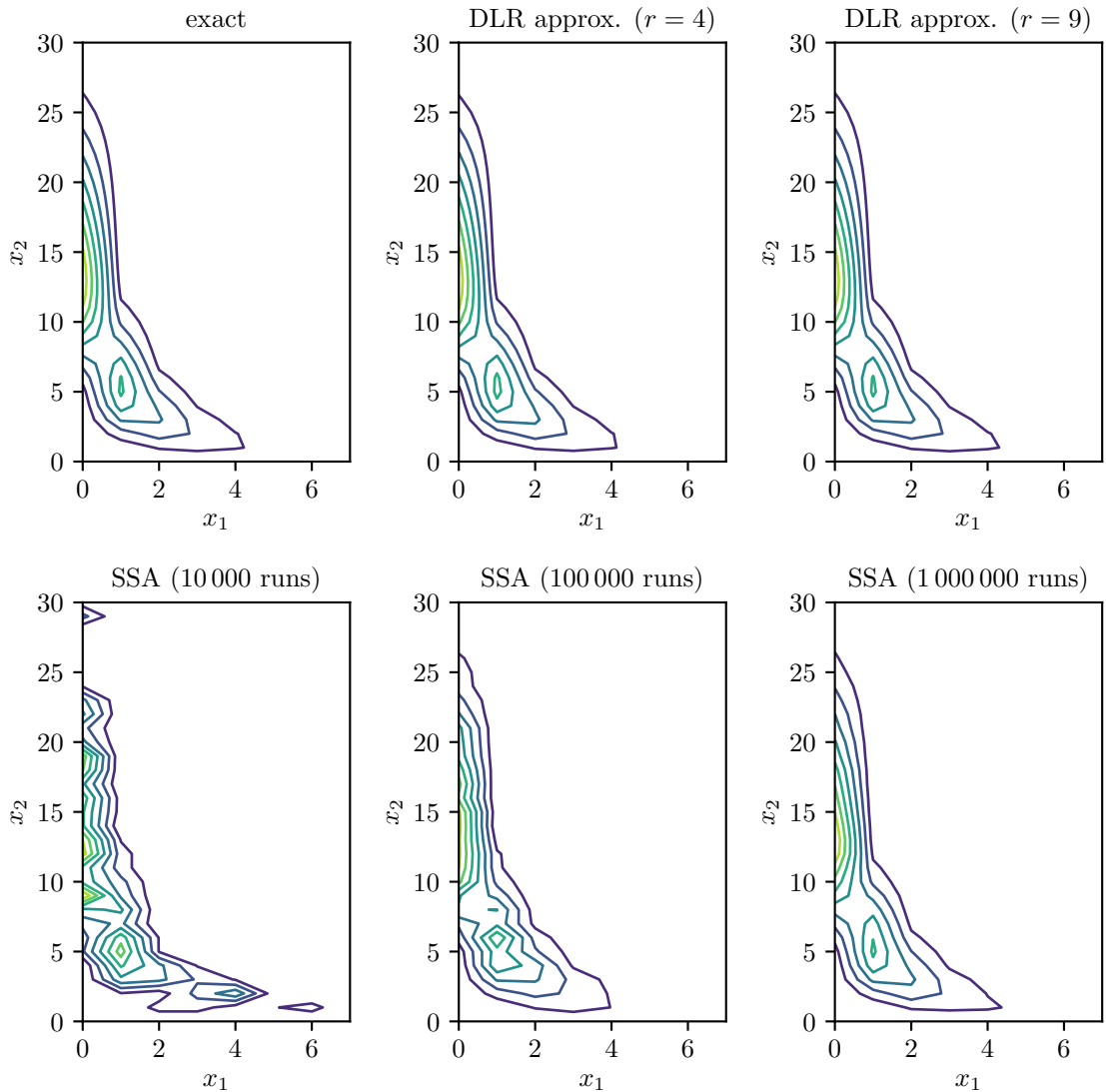}

\caption{Partially evaluated two-dimensional probability distribution $P_{\mathrm{S}}(x_{1},x_{2})=P(x_{1},x_{2},x_{3}=1,x_{4}=1,x_{5}=1)$
of the lambda phage example at $t=10$. Top row: exact reference solution
(solution of the full CME on the truncated state space with the \texttt{scipy.solve\_ivp}
routine) and solutions of the DLR approximation with rank $r=4$ and
$r=9$. The DLR solutions were computed with the second-order integrator
using time step size $\tau=0.01$ and $10$ substeps. Bottom row:
results obtained with SSA using $10\,000$, $100\,000$ and $1\,000\,000$
runs. Note that $P_{\mathrm{S}}$ is only defined at the discrete
grid points $x\in\mathbb{N}_{0}^{2}$; the contour plots are based
on interpolation and are shown here for the sake of clarity.\label{fig:lp2}}
\end{figure}

The stochastic noise of SSA is also visible in figure~\ref{fig:lp2}.
Here the partially evaluated two-dimensional probability distribution
$P_{\mathrm{S}}(x_{1},x_{2})=P(x_{1},x_{2},x_{3}=1,x_{4}=1,x_{5}=1)$
of the lambda phage example is shown for time $t=10$. The probability
distributions calculated with the DLR approximation for rank $r=4$
and $r=9$ agree very well with the exact reference solution, whereas
the noise for the SSA results with $10\,000$ and $100\,000$ samples
is very pronounced. Only for $1\,000\,000$ runs the results are comparable
to the ones obtained by the DLR approximation.

\begin{figure}[H]
\centering
\input{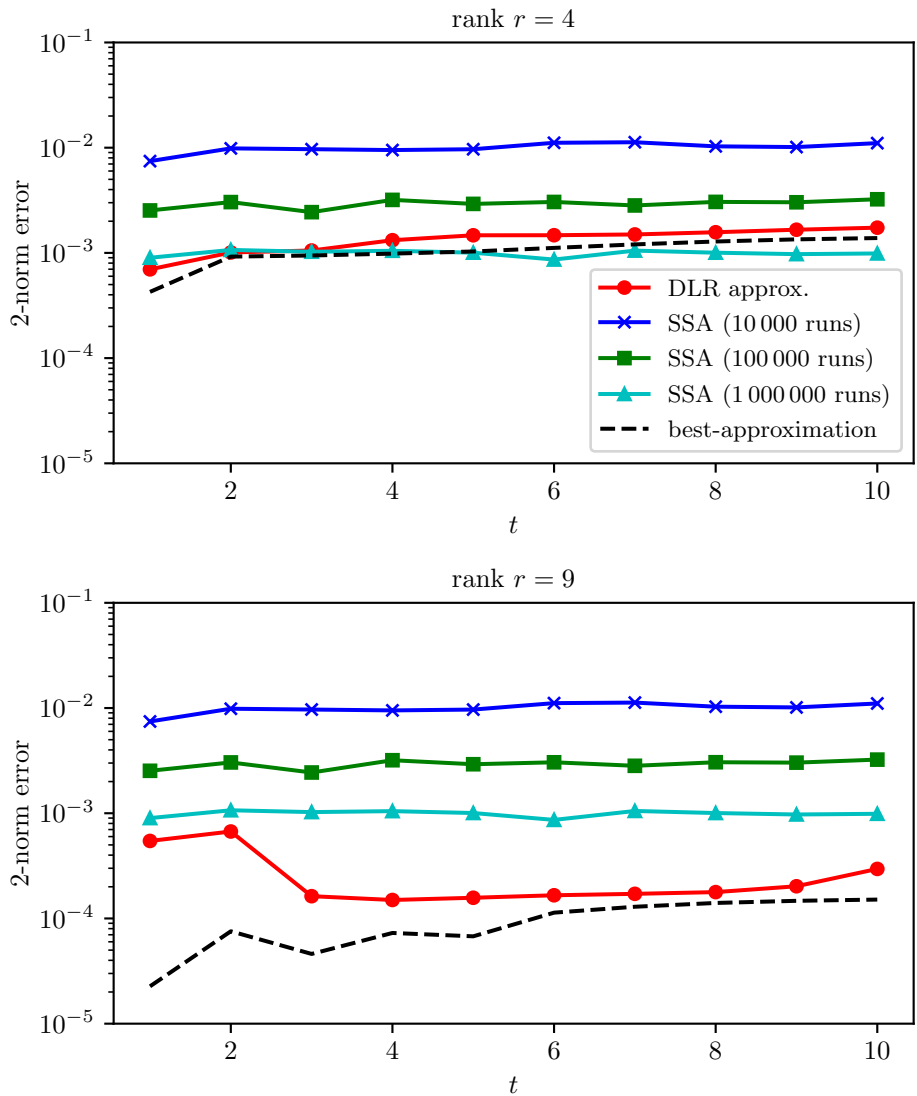}

\caption{Comparison of the 2-norm error as a function of time $t$ for the
DLR approximation (red circles) with the best-approximation (black,
dashed line) using ranks $r=4$ (top) and $r=9$ (bottom) and with
SSA for $1000$ (blue crosses), $100\,000$ (green squares) and $1\,000\,000$
samples (cyan triangles) for the lambda phage example. The error was
calculated by comparing each approximation with the exact reference
solution of the lambda phage, which was obtained by solving the full
CME on the truncated state space directly with the \texttt{scipy.solve\_ivp}
routine. The best-approximation was obtained by truncating all but
the first $r=4$ or $r=9$ singular values for a SVD of the exact
reference solution. For the DLR approximation the second-order integrator
with time step size $\tau=0.01$ and $10$ substeps was used. Note
that the 2-norm error for $t=0$ is not shown in the plot, as it is
zero up to machine precision for the best-approximation, DLR approximation
and SSA with $1\,000\,000$ runs.\label{fig:lp3}}
\end{figure}

Figure~\ref{fig:lp3} shows the 2-norm error of the probability density
function for the DLR approximation, SSA and the best-approximation
depending on time $t$. For rank $r=4$ the resulting error of the
DLR and the best-approximation is slightly larger than for the SSA
using $1\,000\,000$ samples. However, DLR in this configuration is
much faster as has been noted before. If the rank is increased to
$r=9$, the DLR approximation also significantly outperforms SSA in
terms of accuracy. 

\subsection{BAX pore assembly}

The last and most challenging example is the BAX pore assembly, which
is a system with 19 reactions and 11 species. The reactions and propensity
functions of this system are listed in table~\ref{tab:bax}. BAX
plays a key role in mediating mitochondrial outer membrane permeabilization
and is therefore a regulator of programmed cell death (\emph{apoptosis}).
The model was taken from \cite{Gaudet_2012} and is part of the\emph{
extrinsic apoptosis reaction model} (EARM, see, e.g., \cite{Albeck_2008}).
Monomeric BAX ($S_{1}$) can assemble to larger complexes ($S_{2}$--$S_{5}$,
reactions 1--5) and the complexes in turn can dissociate (reactions
6--10). Large enough complexes are able to transport cargo ($S_{10}$
and $S_{11}$), this process is described by reactions 11--19.
\begin{table}[H]
\begin{centering}
\begin{tabular}{c|c|cc}
\hline 
No. & Reaction & Propensity function & \tabularnewline
\hline 
1 & $S_{1}+S_{1}\longrightarrow S_{2}$ & $a_{f}\cdot x_{1}(x_{1}-1)/2$ & \tabularnewline
2--5 & $S_{i}+S_{1}\longrightarrow S_{i+1}$ & $a_{f}\cdot x_{i}x_{1}$ & $(i=2,\dots,5)$\tabularnewline
6--10 & $S_{j+1}\longrightarrow S_{j}+S_{1}$ & $a_{r}\cdot x_{j+1}$ & $(j=1,\dots,5)$\tabularnewline
11--13 & $S_{k}+S_{10}\longrightarrow S_{k+3}$ & $b_{f}\cdot x_{k}x_{10}$ & $(k=4,5,6)$\tabularnewline
14--16 & $S_{k+3}\longrightarrow S_{k}+S_{10}$ & $b_{r}\cdot x_{k+3}$ & \tabularnewline
17--19 & $S_{k+3}\longrightarrow S_{k}+S_{11}$ & $c_{r}\cdot x_{k+3}$ & \tabularnewline
\hline 
\end{tabular}
\par\end{centering}
\caption{Reactions and propensity functions of the BAX pore assembly system
with parameters $a_{f}=2\cdot10^{-4}$, $a_{r}=b_{r}=10^{-3}$, $b_{f}=3\cdot10^{-5}$
and $c_{r}=10$.\label{tab:bax}}
\end{table}

We solved the CME on the time interval $[0,145]$ with rank $r=5$
on the truncated state space with truncation indices $\eta=(0,0,0,0,0,0,0,0,0,0,0)$
and $\zeta=(46,16,16,11,11,11,4,4,4,56,56)$. Note that the purpose
of this numerical example was to discover possible limitations of
our approach; in order to reach the equilibrium one would have to
consider a substantially longer interval of approximately $[0,20\,000]$.
The reaction network was partitioned into $\mathcal{P}_{1}=\{S_{1},S_{2},S_{3},S_{4},S_{5}\}$
and $\mathcal{P}_{2}=\{S_{6},S_{7},S_{8},S_{9},S_{10},S_{11}\}$,
therefore the two partitions have $n_{1}=46\cdot16^{2}\cdot11^{2}=1\,424\,896$
and $n_{2}=11\cdot4^{3}\cdot56^{2}=2\,207\,744$ degrees of freedom.
Thus the total number of degrees of freedom is reduced from $n_{1}\cdot n_{2}=3.15\cdot10^{12}$
to $(n_{1}+n_{2})\cdot r+r^{2}=18\,163\,225$, which is a reduction
by a factor of approximately $1.7\cdot10^{5}$.

We consider the following initial distribution
\[
P(0,x)=\gamma\cdot\exp\left(-\frac{1}{2}(x-\mu)^{T}C^{-1}(x-\mu)\right),
\]
with $C=0.2$, $\mu=(40,0,0,0,0,0,0,0,0,50,0)$ and $\gamma$ was
determined by the condition that $\sum_{x\in\Omega^{\zeta,\eta}}P(0,x)=1$.
We performed the computations for the DLR approximation with the second-order
integrator with $100$ substeps and using a variable time step size.
This time step size was adjusted according to the maximal reaction
rate obtained by solving the rate equations deterministically (which
is very cheap); the minimal time step size is $\tau=1.0$. Due to
the large system size, solving the full CME on the truncated state
space was clearly not possible (we would need approximately $50$
TB of main memory). For comparison we thus consider SSA simulations
with \emph{StochKit2}. The total run time of the DLR approximation
and SSA computations is listed in table \ref{tab:bax_runtime}.
\begin{table}[H]
\begin{centering}
\begin{tabular}{c|r}
\hline 
 & run time {[}s{]}\tabularnewline
\hline 
DLR approx. ($r=5$) & $1.3\cdot10^{5}$\tabularnewline
SSA ($10\,000$ runs) & $84$\tabularnewline
SSA ($100\,000$ runs) & $129$\tabularnewline
SSA ($1\,000\,000$ runs) & $358$\tabularnewline
SSA ($10\,000\,000$ runs) & $2898$\tabularnewline
\hline 
\end{tabular}\caption{Overview of the approximate run times in seconds for the DLR approximation
and SSA for the BAX pore assembly system. The DLR approximation was
computed with the second-order integrator using a variable time step
size and $100$ substeps. All computations were performed on a workstation
with a $2.9$ GHz Intel Core i5 Comet Lake (10400F) processor. The
results of the DLR approximation were computed with six threads.\label{tab:bax_runtime}}
\par\end{centering}
\end{table}

\begin{figure}[H]
\centering
\input{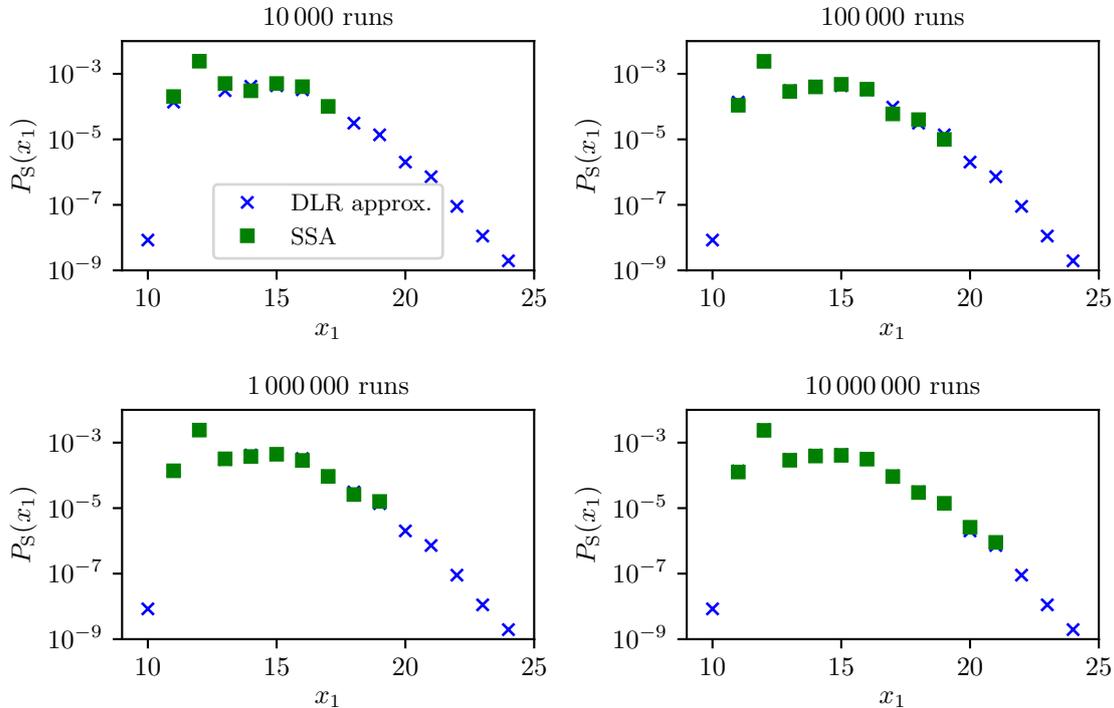}

\caption{Partial evaluated probability distribution $P_{\mathrm{S}}(x_{1})=P(x_{1},x_{2}=9,x_{3}=2,x_{4}=1,x_{5}=0,x_{6}=0,x_{7}=0,x_{8}=0,x_{9}=0,x_{10}=50,x_{11}=0)$
of the BAX pore assembly system at $t=145$. The results were obtained
with our implementation of the DLR approximation using rank $r=5$
(red dots) and with SSA (green squares) using $10\,000$, $100\,000$,
$1\,000\,000$ and $10\,000\,000$ runs. For the DLR approximation
the second-order integrator with variable time step size and $100$
substeps were used.\label{fig:bax1}}
\end{figure}

Figure~\ref{fig:bax1} shows the partially evaluated probability
distribution $P_{\mathrm{S}}(x_{1})=P(x_{1},x_{2}=9,x_{3}=2,x_{4}=1,x_{5}=0,x_{6}=0,x_{7}=0,x_{8}=0,x_{9}=0,x_{10}=50,x_{11}=0)$
at time $t=145$ computed with the DLR approximation using rank $r=5$
and with SSA using $10\,000$, $100\,000$, $1\,000\,000$ and $10\,000\,000$
runs. The results of both methods are in good agreement, which demonstrates
that in principle even such large problems can be solved with our
implementation of the DLR approximation. Even though we use a large
number of samples, SSA has again problems to resolve the tail of the
distribution. Although, we have no exact solution and thus can not
confirm this with certainty, the tail of the dynamical low-rank approximation
follows a power law that looks correct.

\begin{figure}[H]
\centering
\input{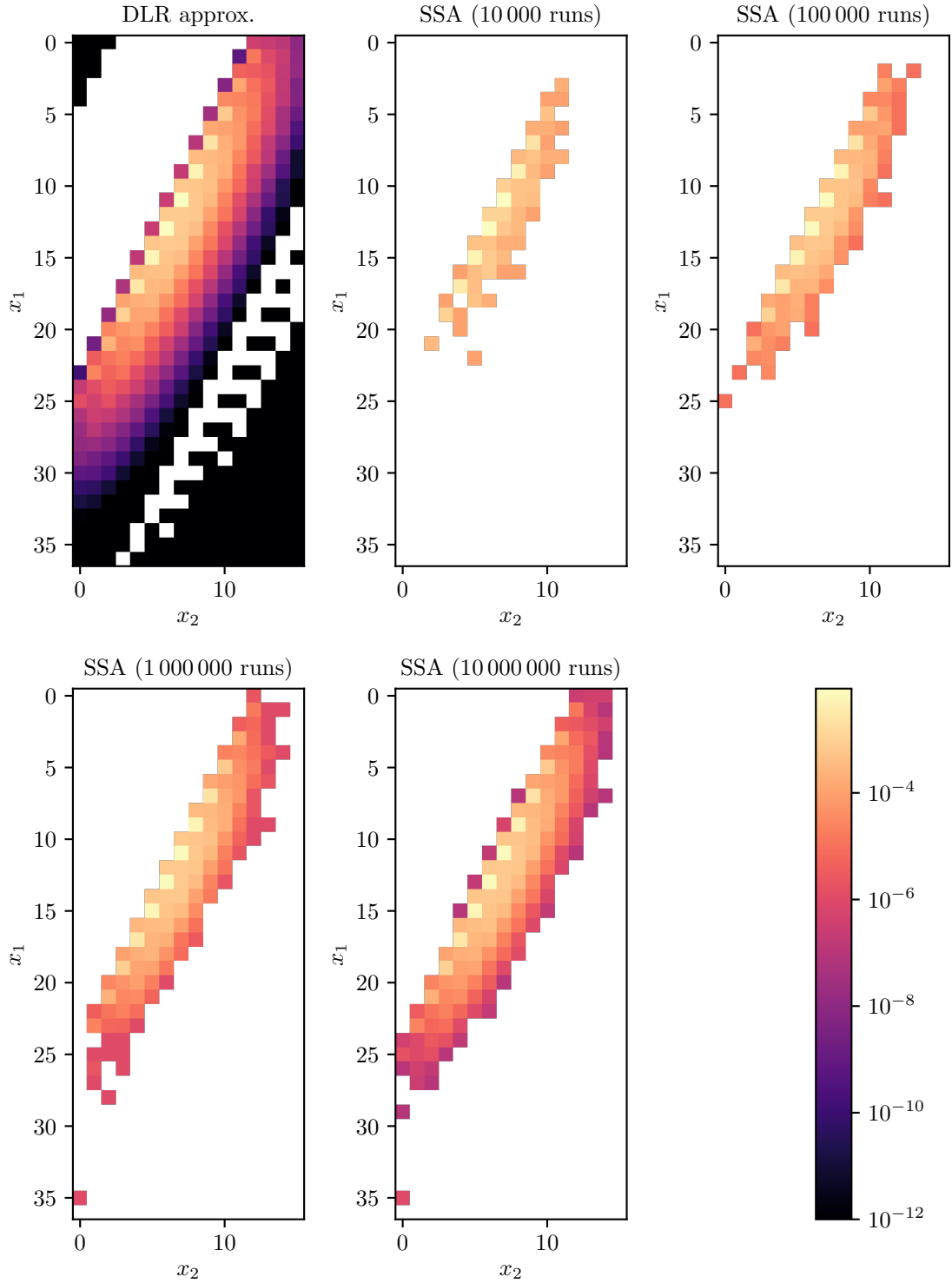}

\caption{Partially evaluated two-dimensional probability distribution $P_{\mathrm{S}}(x_{1},x_{2})=P(x_{1},x_{2},x_{3}=2,x_{4}=1,x_{5}=0,x_{6}=0,x_{7}=0,x_{8}=0,x_{9}=0,x_{10}=50,x_{11}=0)$
of the BAX pore assembly system at $t=145$. The results were obtained
with our implementation of the DLR approximation using rank $r=5$
and with SSA using $10\,000$, $100\,000$, $1\,000\,000$ and $10\,000\,000$
runs. For the DLR approximation the second-order integrator with variable
time step size and $100$ substeps were used. Note that for coloring
a logarithmic mapping was employed; white areas indicate very small
negative (and therefore unphysical) results for the DLR approximation
and zero events in the case of SSA.\label{fig:bax2}}
\end{figure}

In figure~\ref{fig:bax2} the partially evaluated two-dimensional
probability distribution $P_{\mathrm{S}}(x_{1},x_{2})=P(x_{1},x_{2},x_{3}=2,x_{4}=1,x_{5}=0,x_{6}=2,x_{7}=0,x_{8}=0,x_{9}=0,x_{10}=50,x_{11}=0)$
is shown for the same setup. For coloring a logarithmic mapping was
employed; white areas indicate very small negative (and therefore
unphysical) results for the DLR approximation and zero events in the
case of SSA. Again, both methods yield similar results for large probability
values, but it can be clearly seen that the DLR approximation captures
areas of low probability which are not present in the SSA results.

Note that in terms of run time, SSA currently beats the DLR approximation
for this large problem. However, when extending the DLR approach to
a hierarchical scheme (where we divide the reaction network into more
than two partitions) where subproblems have a similar size as the
lambda phage problem, we expect that the run time would be comparable
to the one for the lambda phage example with the additional benefit
that the solutions are noise-free. We consider this the subject of
future work.

\section{Conclusion and outlook\label{sec:Outlook}}

The present work shows that using dynamical low-rank approximations
can result in an algorithm that drastically reduces the memory and
computational effort that is required in order to solve the chemical
master equation. The proposed approach can even outperform SSA by
a significant margin. It is further interesting to note that the DLR
approach directly provides a low-storage approximation of the probability
distribution function (which in SSA has to be reconstructed from the
samples collected as a post-processing step).

The present work considers dividing the problem into two partitions.
However, for large problems this is not sufficient in order to reduce
the memory requirement and computational time to an acceptable level
(and thus to outperform SSA). Thus, as future work, we will consider
the techniques in \cite{Lubich_2013,Einkemmer_2018,Ceruti_2020,Ceruti_2022b}
in order to extend the proposed method to a hierarchical division
into multiple partitions.

One significant advantage of the dynamical low-rank approach considered
here is that it lends itself very well to implicit methods (compared
to, e.g., a step-truncation low-rank approach as considered in \cite{Allmann-Rahn_2022,Cai_2017,Guo_2022,Kormann_2015}).
This is a significant advantage for solving the CME as reaction networks
often include reactions with widely disparate time scales, thus making
the resulting equations stiff. We note that this is also an issue
for SSA (see, e.g., \cite{Harris_2006}).

\bibliographystyle{plain}
\bibliography{bibliography}

\end{document}